\documentclass[12pt]{amsart}
\usepackage{amsmath, amssymb, latexsym, amsthm, mathrsfs, color, mathtools, xcolor, lmodern}
\usepackage[hyphens]{url}
\usepackage[all]{xy}

\usepackage[margin=1in]{geometry}

\newcommand{\nc}{\newcommand}
\newcommand{\vin}{\rotatebox[origin=c]{90}{$\in$}}
\nc{\gm}[1]{ \stackrel{\mbox{\normalsize $#1$}}{\tiny\vin} }
\nc{\slt}{\prec}
\nc{\sqe}{\sqsubseteq^*}
\nc{\splus}{\!+\!}
\nc{\UnN}[1]{\Un_{#1=1}^\oo}
\nc{\CapN}[1]{\bigcap_{#1=1}^\infty}
\nc{\AddCite}{\Cite{Add Citation}}
\nc{\addref}{\Cite{Add Ref}}
\nc{\sfoo}{\sfin(\Op,\Op)}
\nc{\ufog}{\ufin(\Op,\Ga)}
\nc{\utwo}{\mathsf{U}_2}
\nc{\utgg}{\utwo(\Ga,\Ga)}
\nc{\utomg}{\utwo(\Om,\Ga)}
\nc{\uid}{\mathsf{U}_{\op{id}}}
\nc{\uidgg}{\uid(\Ga,\Ga)}
\nc{\sogg}{\sone(\Ga,\allowbreak\Ga)}
\nc{\omg}{\smallbinom{\Om}{\Ga}}

\DeclareMathOperator{\dist}{dist}
\nc{\thusfar}{\par\bigskip\centerline{\textcolor{red}{--- Edited thus far ---}}\par\bigskip}
\nc{\lei}{\le^\oo}
\nc{\card}[1]{\left|#1\right|}
\nc{\medcard}[1]{\biggl|\,#1\,\biggr|}
\nc{\smallcard}[1]{|\,#1\,|}
\nc{\bds}{bidirectional $\roth$-scale}
\nc{\bbN}{\mathbb{N}}
\nc{\mbq}{\mb{?}}
\nc{\mb}[1]{{\mbox{\textbf{#1}}}}
\nc{\nop}{$\times$}
\nc{\fbn}{\!\!\fbox{\!\nop\!}\!\!}
\nc{\yup}{\checkmark}
\nc{\forces}{\Vdash}
\nc{\name}[1]{\dot{#1}}
\nc{\tf}{\my{FINISHED THUS FAR}}
\nc{\FU}{Fr\'echet--Urysohn}
\nc{\gs}{$\gamma$~space}
\nc{\Ga}{\Gamma}\nc{\Om}{\Omega}
\nc{\smallbinom}[2]{\begin{psmallmatrix} #1\\ #2 \end{psmallmatrix}}
\nc{\bgamma}{\smallbinom{\Om}{\Ga}}
\newcommand{\two}{\{0,1\}}
\nc{\productive}[2]{\bigl(#1,\allowbreak #2\bigr)^\x}
\nc{\Sel}{\mathsf{S}}
\nc{\sset}[2]{\{\,#1 : #2\,\}}
\nc{\smb}[1]{{\!\!\mb{#1}\!\!}}
\nc{\medset}[2]{{\biggl\{\,#1 : #2\,\biggr\}}}
\nc{\smallmedset}[2]{{\bigl\{\,#1 : #2\,\bigr\}}}
\nc{\set}[2]{{\left\{\,#1 : #2\,\right\}}}
\nc{\seq}[1]{#1_1,#1_2,\dotsc}
\nc{\eseq}[1]{#1_1, \allowbreak #1_2, \allowbreak\dots} 
\nc{\cube}{(\Cantor)^\bbN}
\nc{\Match}{\op{Match}}
\nc{\concat}[1]{\hat{\phantom{a}}\langle #1\rangle}
\nc{\poset}{\mathbb{P}}
\nc{\fn}[1]{{\op{Fn}(#1\times\w,2)}}
\nc{\linadd}{\op{linadd}}
\nc{\nonprod}{\non^\x}
\nc{\alephes}{{\aleph_0}}
\nc{\my}[1]{\marginpar{\textcolor{red}{***}}\textcolor{red}{#1}}
\nc{\later}[1]{{\color{green} #1}}
\nc{\BTs}[1]{{\color{green} #1 (BT)}}
\nc{\Cp}{\op{C}}
\nc{\Bp}{\op{B}_p}
\nc{\Pa}[8]{\bibitem{#1} {#2}, \emph{#3}, {#4} \textbf{#5} ({#6}), {#7}--{#8}.}
\nc{\tPa}[5]{\bibitem{#1} {#2}, \emph{#3}, {#4}, to appear.}
\nc{\sPa}[4]{\bibitem{#1} {#2}, \emph{#3}, {#4}, submitted.}
\nc{\Bc}[9]{\bibitem{#1} {#2}, \emph{#3}, in: \textbf{#4} (#5), #6 #7, #8--#9.}
\nc{\fD}{\mathfrak{D}}
\nc{\fX}{\mathfrak{X}}
\nc{\Onbd}{\Op_{\mathrm{nbd}}} 
\nc{\Omnb}{\Om_{\mathrm{nbd}}} 
\nc{\od}{\mathfrak{od}}
\nc{\Setting}[7]{\xymatrix@R=4pt@C=7pt{#1\ar@{-}[r]&#2\ar@{-}[r]&#3\\&#4\ar@{-}[u]\\
#5\ar@{-}[uu]\ar@{-}[r] & #6\ar@{-}[u]\ar@{-}[r] & #7\ar@{-}[uu]}}
\nc{\mx}[1]{\begin{matrix}#1\end{matrix}}
\nc{\plim}{p\txt{-}\lim}
\nc{\Bgp}{{\Z^\bbN}}
\nc{\Cgp}{{{\Z_2}^\bbN}}
\nc{\Cite}[1]{\textcolor{red}{[#1]}}
\nc{\Next}[1]{{#1^+}}
\nc{\cFin}{\mathrm{cF}}
\nc{\obbl}{\boldsymbol{\left(\right.}}
\nc{\obbr}{\boldsymbol{\left.\right)}}
\nc{\cbbl}{\boldsymbol{\left[\right.}}
\nc{\cbbr}{\boldsymbol{\left.\right]}}
\nc{\ooi}[2]{\obbl #1(#2),\allowbreak #1(#2\splus 1) \obbr}
\nc{\coi}[2]{\cbbl #1(#2),\allowbreak #1(#2\splus 1) \obbr}
\nc{\cci}[2]{\cbbl #1(#2),\allowbreak #1(#2\splus 1) \cbbr}
\nc{\ooix}[3]{\obbl #1(#2),\allowbreak #1(#3) \obbr}
\nc{\Bdd}{\mathbf{B}}
\nc{\Dfin}{\mathfrak{D}_\mathrm{fin}}
\nc{\grbl}{{\mbox{\textit{\tiny gp}}}}
\nc{\bbP}{\mathbb{P}}
\nc{\BOfat}{\B_{\Om_{\mathrm{fat}}}}
\nc{\Bgood}{\B_{\mathrm{good}}}
\nc{\compactN}{\cl{\mathbb{N}}}
\nc{\blocks}[2]{\op{cl}_{#2}(#1)}
\nc{\blocksplus}[2]{\op{cl}^+_{#2}(#1)}
\nc{\arx}[1]{\texttt{http://arxiv.org/math/#1}}
\nc{\cl}[1]{\overline{#1}}
\nc{\CH}{the Continuum Hypothesis}
\nc{\MA}{Martin's Axiom}
\nc{\Bfat}{\B_\mathrm{fat}}
\nc{\inv}{^{-1}}
\nc{\Cantor}{{\two^\bbN}}
\nc{\bP}{\mathbf{P}}
\nc{\bof}{\op{\fb}}
\nc{\dof}{\op{\fd}}
\nc{\bofF}{\bof(\cF)}
\nc{\sr}[3]{\underset{\mbox{#3}}{\mbox{#1}}}
\nc{\gp}{\binom{\Om}{\Ga}}
\nc{\gpsmall}{\mbox{$\gp$}}
\nc{\gig}{\gimel}
\nc{\gns}{\sone(\Om,\gig)}
\nc{\nsr}[2]{#1}
\nc{\Srg}{{\mathbb{S}}}
\nc{\Srgs}{{\mathbb{S}^*}}
\nc{\NN}{{\bbN^{\bbN}}}
\nc{\ZN}{{\Z^{\bbN}}}
\nc{\NNup}{{\bbN^{\uparrow\bbN}}}
\nc{\Pof}{\op{P}}
\nc{\PN}{{\Pof(\bbN)}}
\nc{\roth}{[\bbN]^{\mbox{\tiny $\infty$}}} 
\nc{\roths}{\textcolor{red}{[\bbN\!\sm\!\{1\}]^{\mbox{\tiny $\infty$}}}}
\nc{\Fin}{\mathrm{Fin}}
\nc{\ici}{[\bbN]^{ \infty, \infty}}
\nc{\Inc}{{\compactN^{\uparrow\bbN}}}
\nc{\powInc}[1]{{\big(\Inc\big)^{#1}}}
\nc{\powFin}[1]{{\big(\Fin\big)^{#1}}}
\nc{\powPN}[1]{{\big(\PN\big)^{#1}}}
\nc{\NcompactN}{{\compactN^\bbN}}
\nc{\Uarrow}{\smash{\big\uparrow}}
\nc{\LE}{\preccurlyeq}
\nc{\GE}{\succcurlyeq}
\nc{\op}{\operatorname}
\nc{\im}{\op{im}}
\nc{\Span}{\op{span}}
\nc{\maxfin}{\op{maxfin}}
\nc{\ran}{\op{range}}
\nc{\iso}{\cong}
\nc{\Madd}{{\M}^\star}
\nc{\cI}{\mathcal{I}}
\nc{\cJ}{\mathcal{J}}
\nc{\scrA}{\mathscr{A}}
\nc{\scrB}{\mathscr{B}}
\nc{\scrC}{\mathscr{C}}
\nc{\scrD}{\mathscr{D}}
\nc{\scrF}{\mathscr{F}}
\nc{\scrK}{\mathscr{K}}
\nc{\A}{\forall}
\nc{\B}{\mathrm{B}}
\nc{\cB}{\mathcal{B}}
\nc{\bB}{\mathbf{B}}
\nc{\BS}{\mathbf{B}(\mathcal{S})}
\nc{\BF}{\mathbf{B}(\mathcal{F})}
\nc{\BU}{\mathbf{B}(\mathcal{U})}
\nc{\cSp}{\mathcal{S}^+}
\nc{\cFp}{\mathcal{F}^+}
\nc{\cUp}{\mathcal{U}^+}

\nc{\BG}{\B_\Ga}
\nc{\BL}{\B_\Lambda}
\nc{\BT}{\B_\Tau}
\nc{\BTstar}{\B_{\Tau^*}}
\nc{\BO}{\B_\Om}
\nc{\DO}{\cD_\Om}
\nc{\KO}{\cK_\Om}
\nc{\CG}{C_\Ga}
\nc{\CL}{C_\Lambda}
\nc{\CT}{C_\Tau}
\nc{\CTstar}{C_{\Tau^*}}
\nc{\CO}{C_\Om}
\nc{\COgp}{C_{\Om^{\grbl}}}
\nc{\CLgp}{C_{\Lambda^{\grbl}}}
\nc{\BOgp}{\B_{\Om}^{\grbl}}
\nc{\BLgp}{\B_{\Lambda^{\grbl}}}
\nc{\sfC}{\mathsf{C}}
\nc{\sfD}{\mathsf{D}}
\nc{\bD}{\mathbf{D}}
\nc{\Tau}{\mathrm{T}}
\nc{\cA}{\mathcal{A}}
\nc{\cK}{\mathcal{K}}
\nc{\cD}{\mathcal{D}}
\nc{\cF}{\mathcal{F}}
\nc{\cS}{\mathcal{S}}
\nc{\cT}{\mathcal{T}}
\nc{\cG}{\mathcal{G}}
\nc{\cY}{\mathcal{Y}}
\nc{\J}{\mathcal{J}}
\nc{\cL}{\mathcal{L}}
\nc{\cM}{\mathcal{M}}
\nc{\cN}{\mathcal{N}}
\nc{\cH}{\mathcal{H}}
\nc{\cO}{\mathcal{O}}
\nc{\Op}{\mathrm{O}}
\nc{\rmA}{\mathrm{A}}
\nc{\rmF}{\mathrm{F}}
\nc{\rmB}{\mathrm{B}}
\nc{\rmD}{\mathrm{D}}
\nc{\rmP}{\mathrm{P}}
\nc{\cC}{\mathcal{C}}
\nc{\cP}{\mathcal{P}}
\nc{\bbQ}{\mathbb{Q}}
\nc{\bbR}{\mathbb{R}}
\nc{\cU}{\mathcal{U}}
\nc{\Un}{\bigcup}
\nc{\cV}{\mathcal{V}}
\nc{\cW}{\mathcal{W}}
\nc{\Z}{{\mathbb Z}}
\nc{\Impl}{\Rightarrow}
\long\def\forget#1\forgotten{}
\nc{\ft}{\mathfrak{t}}
\nc{\fb}{\mathfrak{b}}
\nc{\fc}{\mathfrak{c}}
\nc{\fd}{\mathfrak{d}}
\nc{\fg}{\mathfrak{g}}
\nc{\oo}{\infty}
\nc{\fr}{\mathfrak{r}}
\nc{\fk}{\mathfrak{k}}
\nc{\bidi}{\mathfrak{bidi}}
\nc{\fu}{\mathfrak{u}}
\nc{\fh}{\mathfrak{h}}
\nc{\fp}{\mathfrak{p}}
\nc{\fj}{\mathfrak{j}}
\nc{\fs}{\mathfrak{s}}
\nc{\w}{\omega}
\nc{\x}{\times}
\nc{\Iff}{\Leftrightarrow}
\newcommand\comp{^{\text{\tt c}}}
\nc{\nin}{\notin}
\nc{\cat}{\hat{\ }}
\nc{\sub}{\subseteq}
\nc{\spst}{\supseteq}
\nc{\sm}{\setminus}
\nc{\subs}{\subseteq^*}
\nc{\sups}{\supseteq^*}
\nc{\les}{\le^*}
\nc{\lesp}{\les_+}
\nc{\lti}{<^{\infty}}
\nc{\leS}{\le_S}
\nc{\leF}{\le_{\mathcal{F}}}
\nc{\leU}{\le_{\mathcal{U}}}
\nc{\rest}{\restriction}

\nc{\la}{\langle}
\nc{\ra}{\rangle}
\nc{\E}{\exists}
\nc{\dom}{\op{dom}}
\nc{\cov}{\op{cov}}
\nc{\add}{\op{add}}
\nc{\cof}{\op{cof}}
\nc{\cf}{\op{cf}}
\nc{\non}{\op{non}}
\nc{\unif}{\op{non}}
\nc{\COV}{\op{COV}}
\nc{\ADD}{\op{ADD}}
\nc{\COF}{\op{COF}}
\nc{\NON}{\op{NON}}
\nc{\impl}{\to}
\nc{\Lp}{\mathcal{L_\p}}
\nc{\Wlog}{without loss of generality}
\newtheorem{thm}{Theorem}[section]
\nc{\bthm}{\begin{thm}} \nc{\ethm}{\end{thm}}
\newtheorem{prop}[thm]{Proposition}
\nc{\bprp}{\begin{prop}} \nc{\eprp}{\end{prop}}
\newtheorem{fact}[thm]{Fact}
\nc{\bfct}{\begin{fact}} \nc{\efct}{\end{fact}}
\newtheorem{prob}[thm]{Problem}
\nc{\bprb}{\begin{prob}} \nc{\eprb}{\end{prob}}
\newtheorem{qtn}[thm]{Question}
\nc{\bqtn}{\begin{qtn}} \nc{\eqtn}{\end{qtn}}
\newtheorem{lem}[thm]{Lemma}
\nc{\blem}{\begin{lem}} \nc{\elem}{\end{lem}}
\newtheorem{claim}[thm]{Claim}
\nc{\bclm}{\begin{claim}} \nc{\eclm}{\end{claim}}
\newtheorem{exs}[thm]{Exercise}
\nc{\bexs}{\begin{exs}} \nc{\eexs}{\end{exs}}
\newtheorem{cor}[thm]{Corollary}
\nc{\bcor}{\begin{cor}} \nc{\ecor}{\end{cor}}
\newtheorem{conj}[thm]{Conjecture}
\nc{\bcnj}{\begin{conj}} \nc{\ecnj}{\end{conj}}
\theoremstyle{definition}
\newtheorem{defn}[thm]{Definition}
\nc{\bdfn}{\begin{defn}} \nc{\edfn}{\end{defn}}
\newtheorem{cnv}[thm]{Convention}
\nc{\bcnv}{\begin{cnv}} \nc{\ecnv}{\end{cnv}}
\newtheorem{obs}[thm]{Observation}
\nc{\bobs}{\begin{obs}} \nc{\eobs}{\end{obs}}
\theoremstyle{remark}
\newtheorem{rem}[thm]{Remark}
\nc{\brem}{\begin{rem}} \nc{\erem}{\end{rem}}
\newtheorem{exam}[thm]{Example}
\nc{\bexm}{\begin{exam}} \nc{\eexm}{\end{exam}}
\nc{\bpf}{\begin{proof}} \nc{\epf}{\end{proof}}
\nc{\be}{\begin{enumerate}}
\nc{\ee}{\end{enumerate}}
\nc{\bi}{\begin{itemize}}
\nc{\bimy}{\my{\begin{itemize}}
\nc{\eimy}{\end{itemize}}}
\nc{\itm}{\item}
\nc{\ei}{\end{itemize}}
\nc{\Subsection}[1]{\goodbreak\subsection*{#1}}

\nc{\sone}{\mathsf{S}_1}
\nc{\sfin}{\mathsf{S}_\mathrm{fin}}
\nc{\ufin}{\mathsf{U}_\mathrm{fin}}
\nc{\Split}{\mathsf{Split}}

\nc{\gone}{\mathsf{G}_1}    \nc{\gfin}{\mathsf{G}_\mathrm{fin}}

\nc{\ed}{

\end{document}
}

\title[Omission of Intervals]{Omission of Intervals:\\[0.5cm]
Deducing covering properties of subsets of the real line from their combinatorial structure}

\newif\ifauthor\authortrue
\ifauthor
\author[B. Tsaban]{Boaz Tsaban}
\address{Boaz Tsaban,
Department of Mathematics, Bar-Ilan University, Ramat Gan, Israel}
\email{tsaban@math.biu.ac.il}
\fi

\begin{document}

\begin{abstract}
We develop a method that we call \emph{omission of intervals},
for establishing topological properties of subsets of the real line based on their combinatorial structure.
Using this method, we obtain conceptual proofs of the  fundamental theorems in this realm, and new results that were hitherto inaccessible.
\end{abstract}

\subjclass{Primary: 54D20; 
Secondary: 03E17. 
}

\keywords{Menger property, Hurewicz property, concentrated set, $\gamma$-set,
scale, tower, Omission of Intervals, selection principles.}

\maketitle

\hfill\emph{This interval is intentionally left blank}


\section{Preface}

We develop \emph{omission of intervals},
a novel method for establishing that sets of real numbers with
appropriate combinatorial structures have strong topological covering properties.
The novelty is not only in the obtained results, but also in the uniform
method used to establish them, and in the flexibility of the proof that makes it possible to easily deduce additional results.
This will be demonstrated thoroughly.

This method is applied in the context of \emph{selection principles}, a program
unifying earlier, separate notions and studies originating from
dimension theory (Menger and Hurewicz), measure theory (Borel), and function spaces (Gerlits--Nagy
and Arhangel'ski\u{\i}).
There are several general introductions to selection principles and their applications~\cite{LecceSurvey, KocSurv, ict, SakaiScheepersPIT, Wiki}.
We will provide the necessary definitions and some historical background as we proceed.
To make the reading easier, we usually do not attribute results that are folklore or immediate.

\ifauthor
I developed parts of this method during 2016 and 2023 courses at Bar-Ilan University and
a minicourse at Sichuan University in Chengdu (2016).
Some pieces of this puzzle were conceived earlier, and presented in my plenary lectures at
the Summer Topology Conference in Leicester (2016),
the European Set Theory Conference (2019),
and TOPOSYM (2022).
I thank the organizers for the invitations and the audiences for their enthusiastic feedback,
Lyubomyr Zdomskyy for his thorough reading of this paper and for his useful comments,
and the referees who generously contributed their time and expertise to review this manuscript.
I am grateful to my wife for her encouragement and support.
\fi

\section{Sets of real numbers}

While theorems about relations among properties are striving for minimal restrictions on the involved topological spaces,
examples of spaces with prescribed properties are best when the spaces are as simple and concrete as possible.
All of our examples will be subsets of the real line.
Since the properties considered here are preserved by continuous images,
we can work in topological spaces homeomorphic to subsets of the real line,
such as the Cantor space and the Baire space.

\subsection{The Cantor space}
Our starting point is Cantor's middle-thirds subset of the real line.
This classic set is constructed as follows.
Let $C_0 := [0,1]$, the closed unit interval.
For each natural number $n$, let
$C_{n}$ be the set obtained by removing the middle third (an open interval) from each interval in the set $C_{n-1}$.
Cantor's set is the intersection $C:=\bigcap_{n=0}^\oo C_n$.

\bdfn
Let $\fc:=\card{\bbR}=2^\alephes$, the cardinality of the continuum.
\edfn

The cardinality of Cantor's set is $\fc$: For each $n$, enumerate the $2^n$ intervals in the set $C_n$ as
\[
\sset{I_{(i_1,\dotsc,i_n)}}{i_1,\dotsc,i_n\in\{0,1\}}.
\]
Then the map
\begin{align*}
\Cantor & \longrightarrow C\\
f & \longmapsto x_f\in\CapN{n} I_{(f(1),\dotsc,f(n))}
\end{align*}
is bijective.

The construction of Cantor's set leads naturally to the following classic fact.

\bprp
\label{prp:cantor}
Every uncountable compact set of real numbers has cardinality $\fc$.
Moreover, if $K$ is an uncountable compact set of real numbers and $D$ is a countable set, then there is an uncountable compact set $J\sub K$ that is disjoint from the set $D$.
\eprp
\bpf
The key to the proof is the following lemma.

\blem
Let $X$ be an uncountable subset of $\bbR$.
There are  disjoint closed intervals $I_1$ and $I_2$ whose intersection with the set $X$ is uncountable.
\elem
\bpf
Fix an interval $[a,b]$ such that the set $[a,b]\cap X$ is uncountable.
Let $a'$ be the supremum of the real numbers $x\ge a$ with
$X\cap[a,x]$ countable, and $b'$ be the infimum of
the real numbers $x\le b$ with $X\cap [x,b]$ countable.
Then $a'<b'$.
Take any intervals $[a',c_1]$ and $[c_2,b']$ with $a'<c_1<c_2<b'$.
\epf

Repeat the construction of Cantor's set, as follows. Let $n$ be a natural number. On stage $n$, we have a set $C_n$ comprising of $2^n$ disjoint closed intervals, each having an uncountable intersection with the compact set $K$.
For each of these intervals $I$, the set $I\cap K$ is an uncountable compact set, and the lemma applies to obtain two
disjoint closed subintervals whose intersection with the set $K$ is uncountable. This gives us $2^{n+1}$ disjoint intervals for the next stage. The injective map is defined as before,
taking intersections of the form $\CapN{n} I_{(f(1),\dotsc,f(n))}\cap K$.

To avoid a countable set $D=\{\seq{x}\}$, notice that on
stage $n$, we can apply the lemma twice to obtain \emph{three} disjoint intervals. Out of these, we pick two that do not contain the point $x_n$, and proceed as above.
\epf

\emph{Cantor's space} is the set $\Cantor$, equipped with the metric
\[
\dist(a,b):=\frac{1}{{\min\sset{n\in\bbN}{a(n)\neq b(n)}}}
\]
for distinct points $a,b\in \Cantor$.
Equivalently, as a topological space, the space $\Cantor$ is equipped with the Tychonoff
product topology.
Accordingly, we have $\dist(a,b)<1/n$ if and only if
\[
	(a(1),\dotsc,a(n))=(b(1),\dotsc,b(n)).
\]
It follows that the complement of an open ball in the Cantor space is a finite union of open balls (of the same radius). Thus, open balls are clopen in the Cantor space.
Moreover, the eventually zero sequences form a countable dense set in the Cantor space,
and thus this space has a countable basis consisting of clopen sets.

The standard bijection from the set $\Cantor$ to the set $C$, introduced above, is a topological homeomorphism.
In particular, the Cantor space $\Cantor$ is homeomorphic to the Cantor set $C$ and,
from a topological perspective,
\emph{every subset of the Cantor space is a set of real numbers}.

\subsection{The Cantor space a the power set}
We identify, via characteristic functions, the set $\Cantor$ with the set $\PN$ of all sets of natural numbers. Since we view the space $\PN$ as a set of real numbers, we denote its points by lowercase letters.
The correspondence between the sets $\Cantor$ and $\PN$ dictates that
the distance between distinct points $a,b\in\PN$ is
\[
\dist(a,b)=\frac{1}{\min a\Delta b}.
\]
Thus, we have $\dist(a,b)<1/n$ if and only if
\begin{equation}
\label{eqn:basic2}
a\cap \{1,\dotsc, n\} = b\cap \{1,\dotsc, n\}.
\end{equation}
There is a natural partition of the Cantor space,
\[
\PN=\roth\cup\Fin,
\]
where $\roth$ is the family of infinite subsets of $\bbN$ and $\Fin$ is the family of finite subsets of $\bbN$.
By Equation~\eqref{eqn:basic2}, or by the correspondence with the set $\Cantor$, the countable set $\Fin$ is dense in the space $\PN$.

\subsection{The Baire space}
The \emph{Baire space} $\NN$ is equipped with the Tychonoff product topology or, equivalently, with the topology generated by the metric
\[
\dist(a,b)=\frac{1}{\min\sset{n\in\bbN}{a(n)\neq b(n)}}.
\]
Thus, the Cantor space $\Cantor$ is a metric subspace of the Baire space.
On the other hand, the Baire space is a \emph{topological} subspace of the Cantor space.
To see this, we first notice that the Baire space $\NN$ is isometric to its subspace consisting of
the increasing functions, by mapping every function
$f$ to the increasing function $g$ defined by
\[
g(n):=f(1)+\dotsb+f(n)
\]
for all $n$.
The \emph{increasing} functions $f\in\NN$ are in bijective correspondence
with their image sets $\im(f)$, an element of the Cantor space $\PN$.
In other words, we identify
every element of the set $\roth$ with its increasing enumeration.
The topology of the space $\roth$, inherited from the Cantor space $\PN$, coincides with the topology induced by its identification with the increasing functions in the Baire space.
The identification described here is not uniformly continuous and,
in particular, not an isometry.
However, it is a homeomorphism and, since we deal here with topological properties, every subset of the Baire space may be viewed as a set of real numbers, too.

We view the elements of the set $\roth$ both as infinite sets and as (increasing) functions,
and we refer to them according to their role in every context. Sometimes, we use
the set $\roth$, on occasions where the set $\NN$ would be equally good.
But, usually, we will use the set $\roth$ when we wish to restrict attention to increasing functions.

\section{Menger's Conjecture}
\label{sec:mc}

\subsection{Menger's property}
Menger's property is one of the oldest, and most important, selective covering properties.
Motivated by studies in dimension theory, Menger defined a certain metric property~\cite{Menger24}.
Soon afterwards, Hurewicz realized that Menger's property is equivalent to the topological covering property provided below~\cite{Hure25}.

To avoid technical trivialities, we consider only nontrivial covers of nontrivial topological spaces.

\bcnv
\label{cnv:cover}
We consider only infinite Hausdorff topological spaces.
Let $X$ be a topological space. A \emph{cover} of the space $X$ is a family $\cU\sub\op{P}(X)$ such that $X=\Un\cU$ but $X\notin\cU$.
More generally, a \emph{cover} of a subset $Y\sub X$ is a family $\cU\sub\op{P}(X)$ such that $Y\sub\Un\cU$, but the set $Y$ is not contained in any single member of the cover.
\ecnv

Recall that a topological space $X$ is \emph{compact} if for every open cover $\cU$ of the space,
there is a finite family $\cF\sub\cU$ that covers the space $X$. Compactness is a useful property,
but even very nice spaces, such as the real line $\bbR$, are not compact.
Menger's property is a natural relaxation of compactness, by considering a \emph{sequence} of open covers instead of just one.

\bdfn
For a topological space $X$, let $\Op(X)$ be the family of open covers of this space.
\edfn

In the definition, the symbol $\Op$ stands for arbitrary open covers.
We will later consider restricted types of covers.

\bdfn
Let $\rmA$ and $\rmB$ be types of covers.
A topological space is $\sfin(\rmA,\rmB)$ if
for each sequence
$\seq{\cU}\in\rmA(X)$, there are finite sets $\cF_1\sub\cU_1$, $\cF_2\sub\cU_2$, \dots,
such that $\UnN{n}\cF_n\in\rmB(X)$.
\emph{Menger's property} is the property $\sfoo$.
\edfn

The selection principle $\sfin$ stands for ``select finite''.
We will consider additional selection principles later.

Every compact space is $\sfoo$.
A topological space is \emph{$\sigma$-compact} if it is a countable union of
compact sets. For a family of sets $\cF$, we denote by $\Un\cF$
(a union without a subscript) the union of the members of
the family $\cF$.

\bprp[Menger]
\label{lem:sigmacompactMen}
Every $\sigma$-compact space is $\sfoo$.
\eprp
\bpf
Assume that $X=\UnN{n}K_n$, a countable union of compact sets.
Given a sequence $\seq{\cU}\in\Op(X)$, choose finite sets
$\cF_1\sub\cU_1$, $\cF_2\sub\cU_2$, \dots, such that
$K_1\sub\Un\cF_1$, $K_2\sub\Un\cF_2$, etc.
\epf

\bcnj[Menger]
A metric space is $\sfoo$ if and only if it is $\sigma$-compact.
\ecnj

We will address this conjecture soon.
The following assertion follows from Hurewicz's formulation of Menger's property.

\blem
\label{lem:Mimages}
\mbox{}
\be
\item Menger's property $\sfoo$ is preserved by continuous images: If a topological space $X$ is $\sfoo$, then for each continuous surjective map $X\to Y$, the space $Y$ is $\sfoo$.
\item Every closed subset of a $\sfoo$ space is $\sfoo$.
\qed
\ee
\elem

\subsection{A combinatorial interpretation of Menger's property}

\bdfn
For functions $a,b\in\NN$, we write $a\le b$ if $a(n)\le b(n)$ for all $n$.
We write $a\les b$ if $a(n)\le b(n)$ for \emph{almost all} $n$, that is,
if the set $\sset{n\in\bbN}{a(n)\le b(n)}$ is cofinite.
We write $a\lei b$ if $a(n)\le b(n)$ for infinitely many $n$.
We define the relations $a <^* b$ and $a<^\oo b$ similarly.
\edfn

Thus, we have $a\nleq^* b$ if and only if $b<^\oo a$.
In practice, it will not matter whether strict or weak inequalities hold,
since we have $a+1\le b$ if and only if $a<b$, and similarly for the other relations.

\bdfn
A set $S\sub\NN$ is \emph{dominating} if
for each function $f\in\NN$ there is a function $s\in S$ with $f\les s$.
\edfn

For example, the set $\roth$ is dominating.

\bdfn
Let $\prec$ be a binary relation on $\NN$.
For a set $S\sub\NN$ and an element
$b\in\NN$, we write $S\prec b$ if $s\prec b$ for all elements $s\in S$.
\edfn

\bcor
Let $S\sub\NN$.
The following assertions are equivalent:
\be
\item The set $S$ is not dominating.
\item $S<^\oo b$ for some function $b\in\NN$.
\item $S\lei b$ for some function $b\in\NN$.\qed
\ee
\ecor

\blem
\label{lem:Minf}
Let $X$ be an $\sfoo$ space.
For each sequence $\seq{\cU}\in\Op(X)$,
there are finite sets $\cF_1\sub\cU_1$, $\cF_2\sub\cU_2$, \dots
such that, for each point $x\in X$,
we have $x\in \Un\cF_n$ for infinitely many $n$.
\elem
\bpf
Split the sequence of open covers into infinitely many disjoint sequences, and apply the property
$\sfoo$ to each of these sequences.
\epf

\blem
\label{lem:rec0}
No $\sfoo$ subset of the Baire space $\NN$ is dominating.
\elem
\bpf
Fix a  $n$. For each $m$, let
\[
O^{n}_m:=\set{f\in\NN}{f(n)=m}.
\]
Let $\cU_n:=\{\seq{O^{n}}\}$.

Assume that a set $S\sub\NN$ is $\sfoo$.
Let $\cF_1\sub\cU_1$, $\cF_2\sub\cU_2$, \dots be finite sets as provided by
Lemma~\ref{lem:Minf}. By adding finitely many elements to each set, we may
assume that, for an appropriate function $g\in\NN$,
we have
\[
\cF_1=\{O^{1}_1,\dotsc,O^{1}_{g(1)}\},
\cF_2=\{O^{2}_1,\dotsc,O^{2}_{g(2)}\},
\cF_3=\{O^{3}_1,\dotsc,O^{3}_{g(3)}\},\dotsc
\]
Then $S\lei g$, and thus the set $S$ is not dominating.
\epf

A topological space is \emph{Lindel\"of} if every open cover of this space has a countable subcover.
Spaces with a countable basis, in particular sets of real numbers, are Lindel\"of.
In particular, the Baire space is Lindel"of.
Thus, in contrast to compact spaces, Lindel"of spaces need not be $\sfoo$.

\bdfn
A cover $\cV$ \emph{refines} a cover $\cU$ if each element of $\cV$ is contained in some element of $\cU$.
\edfn

\blem
\label{lem:cloprefine}
Let $X$ be a subset of the Cantor space. Every open cover of $X$ is refined by a countable clopen cover whose members are pairwise disjoint.
\elem
\bpf
Recall that the Cantor space has a countable basis whose elements are clopen sets.
Let $\cU\in\Op(X)$. Replacing every open set $U\in\cU$ by all of its basic clopen subsets,
the family $\cU$ becomes a countable cover by clopen sets.
Enumerate this new family, $\cU=\{\seq{C}\}$.
Replacing each clopen set $C_n$ by the clopen set $C_n\sm (C_1\cup\dotsb\cup C_{n-1})$, we may
assume that the elements of the cover $\cU$ are also pairwise disjoint.
\epf

The following theorem is implicit in Hurewicz's paper~\cite{Hure27}.
In the present form, it was established by Rec\l{}aw~\cite[Proposition~3]{Rec94}.

\bthm[Rec\l{}aw]
Let $X$ be a subset of the Cantor space. The following assertions are equivalent:
\be
\item The set $X$ is $\sfoo$.
\item No continuous image of the set $X$ in the Baire space $\NN$ is dominating.
\ee
\ethm
\bpf
$(1)\Impl (2)$: The property $\sfoo$ is preserved by continuous images (Lemma~\ref{lem:Mimages}),
and $\sfoo$ subsets of $\NN$ are not dominating (Lemma~\ref{lem:rec0}).

$(2)\Impl (1)$:
By Lemma~\ref{lem:cloprefine}, it suffices to consider countable disjoint clopen covers.
Assume that, for each $n$,
the family
$\cU_n=\{\seq{U^n}\}$
is a clopen disjoint cover of $X$.
Define a map
\begin{align*}
 X & \longrightarrow  \NN\\
 x & \longmapsto f_x(n):= m\text{ such that }x\in U^n_m.
\end{align*}
This map is continuous: Fix a natural number $n$.
For each point $x\in X$, we have $x\in U^n_{f_x(n)}$.
If a point $y\in X$ is close enough to the point $x$, then $y\in U^n_{f_x(n)}$,
and thus $f_y(n)=f_x(n)$.

By (2), the image of the map $\Psi$ is not dominating, and thus there is a function $g\in\NN$
such that $\Psi(X)\lei g$.
Then $\UnN{n}\{U^n_1,\dotsc,U^n_{g(n)}\}\in\Op(X)$.
\epf

The set $\roth$, viewed as a subspace of the Baire space $\NN$, is dominating.
By Rec\l{}aw's Theorem, dominating subsets of the space $\roth$ are not $\sfoo$.

\bdfn
Let $\fd$ be the minimal cardinality of a dominating set $S\sub\NN$.
\edfn

Since the entire set $\NN$ is dominating, we have $\fd\le\fc$. It is also
easy to see that a countable subset of $\NN$ cannot be dominating.
Thus, $\aleph_1\le\fd\le\fc$.
No equality is provable here; equivalently,
strict inequalities are consistent~\cite{BlassHBK}.
Nonetheless, we will use the cardinal $\fd$ to obtain results in ZFC!

\bdfn
Let $\bP$ be a property of some, but not all, subsets of the Cantor space.
The \emph{critical cardinality} of $\bP$, denoted
$\non(\bP)$, is the minimal cardinality of a subset of the Cantor space
that does not have the property $\bP$.
\edfn

Rec\l{}aw's Theorem implies the following result.

\bcor
\label{cor:nonsfoo}
$\non(\sfoo)=\fd$.
\ecor
\bpf
$(\ge)$ By Rec\l{}aw's Theorem, an example for $\non(\sfoo)$ has a dominating image in the Baire space, so its cardinality is at least $\fd$.

$(\le)$ The space $\NN$ is a subspace of the Cantor space.
Use Lemma~\ref{lem:rec0}.
\epf

\subsection{Addressing Menger's Conjecture}

\bdfn
Let $\kappa$ be an uncountable cardinal.
A set $S\sub\roth$ is \emph{$\kappa$-unbounded} if $\card{S}\ge\kappa$, and for each function $b\in\NN$,
we have $\card{\sset{s\in S}{s\le b}}<\kappa$.
\edfn

\brem
\label{rem:unccof}
A cardinal $\kappa$ has \emph{uncountable cofinality} if it is not the limit of a countable increasing sequence of cardinals.
We will consider $\kappa$-unbounded sets for cardinal numbers $\kappa$ of uncountable cofinality.
For these cardinals, changing $\le$ to $\les$ will not alter our definition.
\erem

If a set $S\sub\roth$ is $\kappa$-unbounded, then this set is $\lambda$-unbounded for all cardinals $\lambda$ with
$\kappa\le\lambda\le\card{S}$.

\blem\label{lem:du}
There are $\fd$-unbounded sets in $\roth$.
\elem
\bpf
Fix a dominating set $\sset{f_\alpha}{\alpha<\fd}\sub\NN$.
For each ordinal $\alpha<\fd$, take an increasing
function $s\in\roth$ with
$\sset{f_\beta}{\beta\le\alpha}<^\oo s$,
and let $s_\alpha$ be an infinite subset of $s$ that is not in the set
$\sset{s_\beta}{\beta<\alpha}$.
Since $s\le s_\alpha$, we have
\[
\sset{f_\beta}{\beta\le\alpha}<^\oo s_\alpha\notin \sset{s_\beta}{\beta<\alpha}.
\]
The set $\sset{s_\alpha}{\alpha<\fd}$ is $\fd$-unbounded:
$\card{S}=\fd$.
For each function $b\in\NN$, there is an ordinal $\alpha$ with $b\les f_\alpha$.
For each ordinal
$\beta\ge\alpha$, we have
\[
b\les f_\alpha <^\oo s_\beta,
\]
and thus $b <^\oo s_\beta$.
\epf

Recall that $\Fin$ is the set of finite subsets of $\bbN$, so that $\PN=\roth\cup\Fin$.

\bdfn
Let $\kappa$ be an uncountable cardinal with $\kappa\le\fc$.
A set $X\sub\PN$ is \emph{$\kappa$-concentrated} on the set $\Fin$
if
$\card{X}\ge\kappa$, and
$\card{X\sm U}<\kappa$ for all open sets $U\spst\Fin$.
\edfn

\blem\label{lem:dconc}
For each set $X\sub\PN$ that is $\fd$-concentrated on $\Fin$,
the set $X\cup\Fin$ is $\sfoo$.
\elem
\bpf
We can cover the countable set $\Fin$ in one pass, and the
remaining, fewer than $\fd$ points, in a second pass. Formally:
Let $\seq{\cU}\in\Op(X)$.
Take finite sets $\cF_1\sub\cU_1$, $\cF_2\sub\cU_2$, \dots
such that $\Fin\sub U:=\Un(\UnN{n}\cF_n)$.
Then $\card{X\sm U}<\fd$, and by Corollary~\ref{cor:nonsfoo}
there are finite sets
$\cF_1'\sub\cU_1$, $\cF_2'\sub\cU_2$, \dots
that together cover the set $X\sm U$.
Then $\UnN{n}(\cF_n\cup\cF_n')\in\Op(X)$.
\epf

\blem\mbox{}
\be
\item For each compact set $K\sub \roth$, there is a function
$b\in\NN$ with $K\le b$.

\item For each function $b\in\NN$, the set $\set{s\in\roth}{s\le b}$ is compact.
\ee
\elem
\bpf
(1) For each $n$, the projection $\sset{f(n)}{f\in K}$ is a continuous image of the compact set
$K$ in $\bbN$. Since a continuous image of a compact set is compact, this projection is finite.

(2) In the Baire space $\NN$, we have
$\sset{s\in\NN}{s\le b}=\prod_n \cbbl 1,b(n)\cbbr$.
By Tychonoff's Theorem, a product of finite sets is compact.
(Alternatively, the set $\sset{s\in\NN}{s\le b}$ is closed and
totally bounded in the complete metric space $\NN$.)

Since the set $\roth$ is closed in the Baire space $\NN$,
the set
\[
\sset{s\in\roth}{s\le b}=\set{s\in\NN}{s\le b}\cap\roth,
\]
a closed subset of a compact set, is compact.
\epf

The following observation of Szewczak and the author~\cite[Lemma~2.3]{pMReal} generalizes one of Lawrence~\cite[Proposition~2]{Lawrence90}.

\bprp
\label{prp:cupfin}
Let $S\sub\roth$. The following assertions are equivalent:
\be
\item The set $S$ is $\kappa$-unbounded.
\item The set $S$ is $\kappa$-concentrated on $\Fin$.
\ee
\eprp
\bpf
$(1)\Impl (2)$: For each open set $U\spst\Fin$,
the set $K:=\PN\sm U\sub\roth$ is compact.
Thus, there is a function $b\in\NN$ such that
$K\le b$.
Then $S\sm U=S\cap K\le b$, and thus
\[
\card{S\sm U}\le\card{\sset{s\in S}{s\le b}}<\kappa.
\]

$(2)\Impl (1)$: Given a function $b\in\NN$, let
$K=\set{s\in\roth}{s\le b}$.
Since the set $K$ is compact, its complement $K\comp$ is open.
As $K\comp\spst \Fin$, we have
\[
\card{\set{s\in S}{s\le b}}=\card{S\cap K}=\card{S\sm K\comp}<\kappa.\qedhere
\]
\epf

\blem\label{lem:nosk}
Assume that $\aleph_1\le\kappa\le\fc$.
If a set $X\sub\PN$ is $\kappa$-concentrated on the set $\Fin$,
then the set $X$ is not $\sigma$-compact.
\elem
\bpf
Let $X\sub\PN$ be $\kappa$-concentrated on $\Fin$.
Since $\card{X}\ge\kappa$, the set $X$ is uncountable.
Assume that it is $\sigma$-compact.
Then some compact set $K\sub X$ is uncountable.
Since the set $\Fin$ is countable, there is an
uncountable compact (and thus of cardinality $\fc$)
set $J\sub K$ disjoint from the set $\Fin$
(Proposition~\ref{prp:cantor}).
The complement $J\comp$ is open and contains the set $\Fin$.
Since the set $X$ is $\kappa$-concentrated on $\Fin$, we have
\[
\fc = \card{J} = \card{X\sm J\comp}<\kappa.
\]
A contradiction.
\epf

Recall that there are $\fd$-unbounded sets in $\roth$ (Lemma~\ref{lem:du}).
We have the following observation~\cite[Remark~18]{ideals}.

\bthm[Bartoszy\'nski--Tsaban]
Let $S\sub\roth$ be a $\fd$-unbounded set.
The set $S\cup\Fin$ is a non-$\sigma$-compact, $\sfoo$ set of real numbers.
\ethm
\bpf
By Proposition~\ref{prp:cupfin}, the set $S\cup\Fin$ is $\fd$-concentrated on $\Fin$.
By Lemma~\ref{lem:dconc}, it is $\sfoo$.
By Lemma~\ref{lem:nosk}, it is not $\sigma$-compact.
\epf

\bdfn
Let $\bP$ be a property of some, but not all, subsets of the Cantor space.
A \emph{nontrivial} $\bP$ set is a set of
cardinality at least $\non(\bP)$ that has the property $\bP$.
\edfn

\bcor
\label{cor:nontM}
There are nontrivial $\sfoo$ sets of real numbers that are not $\sigma$-compact.
\ecor

\subsection{Comments for Section~\ref{sec:mc}}
\label{subsec:mc}
A set of real numbers is \emph{analytic} if it is a continuous image of a Borel set.
Hurewicz~\cite{Hure25} proved that when restricted to analytic sets, Menger's Conjecture
is true.
Indeed, he proved that an analytic set either contains a closed subset homeomorphic to the Baire space, or else it is $\sigma$-compact.
Since Menger's property is hereditary for closed subsets, and the Baire space is not Menger, the assertion follows.

A game-theoretic variation of Menger's Conjecture is also true:
The \emph{Menger game}~\cite{Hure25}, $\gfin(\Op,\Op)$, is the game associated to Menger's property $\sfin(\Op,\Op)$.
It is played on a topological space $X$, and has an inning per each natural number $n$.
In each inning, Alice picks an open cover $\cU_n$ of the space, and Bob chooses a finite set $\cF_n\sub\cU_n$.
Bob wins if the family $\UnN{n}\cF_n$ is a cover of the space, and otherwise Alice wins.
Telg\'arsky~\cite{TelGames3, Telgarsky87} proved that if Bob has a winning strategy in the game $\gfin(\Op,\Op)$ played on a metric space, then the space is $\sigma$-compact.
Scheepers~\cite[Theorem~1]{Sch95} provides a direct proof, that was slightly simplified recently~\cite[Theorem~4.3]{BKPfs}.

A set of real numbers is \emph{Luzin} if its intersection with every nowhere dense set is countable.
The Continuum Hypothesis implies the existence of Luzin sets.
Sierpi\'nski~\cite{Sier26} observed that every Luzin set is a counterexample for Menger's Conjecture.
In our terminology, every Luzin set is $\aleph_1$-concentrated on a countable set, and is thus $\sfoo$ but not $\sigma$-compact.

Fremlin and Miller~\cite[Theorem~4]{FM88} were the first to prove, outright in ZFC, that Menger's Conjecture is false.
Their argument is dichotomic: It provides one counterexample in the case $\aleph_1=\fd$ and another, trivial counterexample in the case $\aleph_1<\fd$.
Corollary~\ref{cor:nontM} makes a stronger assertion.

\section{Interlude: Selecting one element from each family}
\label{sec:RothS1GO}

\subsection{Rothberger's property}

Menger's property is an instance of the selection principle $\sfin$.
We now introduce a related selection principle.

\bdfn
Let $\rmA$ and $\rmB$ be types of covers of topological spaces.
A topological space $X$ is $\sone(\rmA,\rmB)$ if for each sequence of covers
$\seq{\cU}\in\rmA(X)$,
there are sets $U_1\in\cU_1$, $U_2\in\cU_2$, \dots
such that $\{\seq{U}\}\in\rmB(X)$.
\edfn

Rothberger~\cite{Roth38} introduced the property $\sone(\Op,\Op)$ as a topological variation of
Borel's \emph{strong measure zero} property~\cite{Borel19}:
A set $X$ of real numbers has strong measure zero if for each sequence of positive real numbers $\seq{\epsilon}$,
there is a cover by intervals $X\sub\UnN{n}I_n$ such that
for each $n$, the length of the interval $I_n$ is smaller than $\epsilon_n$.
Every $\sone(\Op,\Op)$ set of real numbers has strong measure zero.
These properties are restrictive:
Borel conjectured that only countable sets have strong measure zero,
and his conjecture is consistent~\cite{Laver76}.
We have the following observation.

\blem
\label{lem:a1conc}
For each set $X\sub\PN$ that is $\aleph_1$-concentrated on $\Fin$,
the set $X\cup\Fin$ is $\sone(\Op,\Op)$.
\elem
\bpf
We modify the proof of Lemma~\ref{lem:dconc}:
We first split the given sequence $\seq{\cU}\in\Op(X)$ to two disjoint subsequences.
We then can cover the countable set $\Fin$ by elements from the first subsequence.
Since the set $X$ is $\aleph_1$-concentrated on $\Fin$,
there remain only countably many uncovered elements.
We cover them with elements from the second subsequence.
\epf

Since the property $\sone(\Op,\Op)$ implies $\sfoo$, we have $\non(\sone(\Op,\Op))\le\non(\sfoo)=\fd$.
Thus, for example, if $\fd=\aleph_1$ then we have a nontrivial $\sone(\Op,\Op)$ of real numbers.
All sets that we consider are $\kappa$-concentrated for some cardinal $\kappa\le\fd$, and thus provide
counterexamples for Borel's conjecture if we assume an appropriate portion of \CH.

\subsection{A provably nontrivial variation of Rothberger's property}

There is another way to see that the sets that we consider are not $\sigma$-compact.

\bdfn
An cover $\cU$ of a topological space $X$ is
\emph{point-cofinite} if it is infinite, and each point $x\in X$ is in almost all sets $U\in\cU$,
that is, for each point $x\in X$, the set $\sset{U\in\cU}{x\in U}$ is a cofinite subset of the cover $\cU$.
Let $\Ga(X)$ be the family of all open point-cofinite covers of the space $X$.
\edfn

Point-cofinite covers are those covers that are hereditary for infinite subsets.

\blem
Every infinite subset of a point-cofinite cover of a topological space is a point-cofinite cover of that space.
Conversely, if all infinite subsets of an infinite cover are also covers, then the cover is point-cofinite.
\qed
\elem

We have the following observation~\cite[Theorem~1.2]{coc2}.

\blem[Just--Miller--Scheepers--Szeptycki]
\label{lem:sfgo}
For Lindel\"of spaces, $\sfin(\Ga,\Op) = \sfoo$.
\elem
\bpf
Let $X$ be an $\sfin(\Ga,\Op)$ space.
We need to prove $\sfoo$.

Let $\seq{\cU}\in\Op(X)$.
If any of these covers has a finite subcover, we are done.
Thus, assume that none has a finite subcover.
For each $n$, thin out the cover $\cU_n$ to obtain a countable cover
\[
\cU_n := \sset{U^n_m}{m\in\bbN}.
\]
Then, for each $m$, replace each set $U^n_m$ with the set
\[
V^n_m:=U^n_1\cup U^n_2\cup\dotsb\cup U^n_m.
\]
The family $\sset{V^n_m}{m\in\bbN}$ is a point-cofinite cover of the space $X$.
By the property $\sfin(\Ga,\Op)$, since the sets $V^n_m$ are increasing with the index $m$,
there are natural numbers $\seq{m}$ with
$\sset{V^n_{m_n}}{n\in\bbN}\in\Op(X)$.
It follows that
\[
\UnN{n}\{U^n_1,U^n_2,\dotsc,U^n_{m_n}\}\in\Op(X).\qedhere
\]
\epf

By Lemma~\ref{lem:sfgo}, for Lindel\"of spaces we have
\[
\sone(\Op,\Op) \longrightarrow \sone(\Ga,\Op) \longrightarrow \sfoo.
\]

We arrive at the following corollary~\cite[Theorem~4.5]{coc2}.

\bcor[Just--Miller--Scheepers--Szeptycki]
\label{cor:nonsogo}
$\non(\sone(\Ga,\Op))=\fd$.
\ecor
\bpf
$(\le)$ The property $\sone(\Ga,\Op)$ implies $\sfoo$, whose critical cardinality is $\fd$ (Corollary~\ref{cor:nonsfoo}).

$(\ge)$ Suppose that $\card{X}<\fd$.
Let $\seq{\cU}\in\Ga(X)$.
By moving to infinite subsets, we may assume that the covers $\cU_n$ are countable, and enumerate
$\cU_n=\sset{U^n_m}{m\in\bbN}$
for each $n$.
For each point $x\in X$, define a function $f_x\in\NN$ by
\[
f_x(n) := \min\sset{m}{x\in U^n_m\text{ for all }m\ge n}.
\]
Let $g\in\NN$ be a function with $\sset{f_x}{x\in X}\lei g$.
Then every point $x\in X$ is in $U^n_{g(n)}$ for infinitely many $n$.
\epf

\blem
\label{lem:dconc2}
For each set $X\sub\PN$ that is $\kappa$-concentrated on $\Fin$ for some cardinal $\kappa\le\fd$,
the set $X\cup\Fin$ is $\sone(\Ga,\Op)$.
\elem
\bpf
The proof is similar to that of Lemma~\ref{lem:a1conc}, using Corollary~\ref{cor:nonsogo}.
\epf

Since there are $\fd$-concentrated sets, we have the following observation~\cite[Remark~18]{ideals}.

\bcor[Bartoszy\'nski--Tsaban]
There are nontrivial $\sone(\Ga,\Op)$ sets of real numbers. \qed
\ecor

There is a fundamental difference between the property $\sone(\Ga,\Op)$ and Menger's property $\sfoo$~\cite[Theorem~2.3]{coc2}.

\bprp[Just--Miller--Scheepers--Szeptycki]
\label{prp:nosk}
The Cantor space is not a subspace of any $\sone(\Ga,\Op)$ space.
\eprp
\bpf
Assume otherwise.
The property $\sone(\Ga,\Op)$ is hereditary for closed subsets, and is preserved by continuous images.
Since a compact subspace is a closed subset, it follows that the Cantor space is $\sone(\Ga,\Op)$.
As topological spaces, we have
\[
\Cantor = \{0,1\}^{\bbN\x\bbN} = (\Cantor)^\bbN.
\]
Fix distinct elements $\seq{x}\in\Cantor$.
For each $n$, define a cover
\[
\cU_n := \sset{U^n_m}{m\in\bbN}\in\Ga( (\Cantor)^\bbN )
\]
by
\[
U^n_m := \sset{f\in (\Cantor)^\bbN}{f(n)\neq x_m}
\]
for all $m$.
Suppose that $\{U^n_{m_1}, U^n_{m_2}, \dotsc\}\in\Op(X)$.
Then the fucntion $f\in (\Cantor)^\bbN$ defined by $f(n):=x_{m_n}$ for all $n$
is not in any set $U^n_{m_n}$.
A contradiction.
\epf

It follows that uncountable $\sone(\Ga,\Op)$ sets of real numbers are not $\sigma$-compact.
The sets we consider are all $\sone(\Ga,\Op)$.

\subsection{Comments for section~\ref{sec:RothS1GO}}
Sierpi\'nski observed that every Luzin set (Section~\ref{subsec:mc}) has strong measure zero.
Indeed, Luzin sets are concentrated on every countable dense subset, and are thus $\sone(\Op,\Op)$.

The critical cardinality of the property $\sone(\Op,\Op)$ is well understood:
Let $\cov(\cM)$ be the minimal cardinality of a set $S\sub\NN$ such that
there is no function $g\in\NN$ with $f =^\oo g$ for all $f\in S$.
(We say that $f =^\oo g$ if $f(n)=g(n)$ for infinitely many $n$.)
This cardinal is named $\cov(\cM)$ since it happens to characterize the minimal cardinality of the minimal cover of the real line by meager sets~\cite{BlassHBK}.
A proof similar to that of Corollary~\ref{cor:nonsfoo} establishes the following result.

\bprp
$\non(\sone(\Op,\Op))=\cov(\cM)$.\qed
\eprp

Lemma~\ref{lem:a1conc} then holds for all $\kappa$-concentrated sets, for $\kappa\le\cov(\cM)$.
We will soon introduce a combinatorial cardinal $\fb$ with $\fb\le\fd$, and see that there are $\fb$-concentrated sets.
It follows that the hypothesis $\fb\le\cov(\cM)$ is sufficient to render Borel's Conjecture false.

\section{The Hurewicz Problem}
\label{sec:hp}

In his study of Menger's property $\sfoo$, Hurewicz~\cite{Hure25} introduced the following, stronger property.

\bdfn
Let $\rmA$ and $\rmB$ be types of covers of topological spaces.
A topological space $X$ is $\ufin(\rmA,\rmB)$ if for each sequence of covers
$\seq{\cU}\in\rmA(X)$, none with a finite subcover,
there are finite sets $\cF_1\sub\cU_1$, $\cF_2\sub\cU_2$, \dots
such that $\{\seq{\Un\cF}\}\in\rmB(X)$.
\emph{Hurewicz's property} is the property $\ufog$.
\edfn

In the notation $\ufin(\rmA,\rmB)$, the selection principle $\ufin$ stands for ``union finite''; we
take the \emph{unions} of the selected finite subfamilies.
Since $\sigma$-compact sets are also
countable \emph{increasing} unions of compact sets, the proof of Lemma~\ref{lem:sigmacompactMen}
establishes the following fact.

\bprp[Hurewicz]
Every $\sigma$-compact space is $\ufog$.\qed
\eprp

We thus have the following implications
\[
\sigma\text{-compact} \longrightarrow \ufog \longrightarrow \sfoo.
\]
Menger's Conjecture asserted that these three properties are equivalent.
Motivated by that, Hurewicz formulated the following milder conjecture and natural problem.

\bcnj[Hurewicz]
\label{cnj:Hur}
Every $\ufog$ metric space is $\sigma$-compact.
\ecnj

\bprb[Hurewicz]
Is there an $\sfoo$ set of real numbers that is not $\ufog$?
\eprb

We will address the conjecture in the next section,
and the problem shortly.
The omitted proofs below are similar to those for $\sfoo$, and are left as an exercise.

\blem
The property $\ufog$ is preserved by continuous images, and is hereditary for closed subsets.\qed
\elem

\subsection{A combinatorial interpretation of Hurewicz's property}

\bdfn
A set $S\sub\NN$ is \emph{bounded} (respectively, \emph{unbounded}) if it is
$\les$-bounded (respectively, $\les$-unbounded).
\edfn

Here, we cannot replace $\les$ by $\le$.
The set of constant functions in $\NN$ is unbounded with respect to the relation $\le$, but it is
$\les$-bounded by the identity function.

\blem
Every $\ufog$ subset of the Baire space $\NN$ is bounded.\qed
\elem

\bthm[Rec\l{}aw]
\label{thm:rec2}
Let $X$ be a subset of the Cantor space.
The following assertions are equivalent:
\be
\item The set $X$ is $\ufog$.
\item Every continuous image of the set $X$ in the Baire space $\NN$ is bounded.\qed
\ee
\ethm

\bdfn
Let $\fb$ be the minimal cardinality of an unbounded subset of $\NN$.
\edfn

\bcor
\label{cor:bineq}
$\aleph_1\le \fb\le\fd\le\fc$.\qed
\ecor

Strict inequalities are consistent in Corollary~\ref{cor:bineq}, but we will later use the cardinal $\fb$, too, to obtain results in ZFC.

\bcor
\label{cor:nonh}
$\non(\ufog)=\fb$.\qed
\ecor

\blem
\label{lem:b-unb}
There are $\fb$-unbounded sets in $\roth$.
\elem
\bpf
Fix an unbounded set $\sset{f_\alpha}{\alpha<\fb}\sub\NN$.
For each ordinal $\alpha<\fb$, take a function
$s_\alpha\in\roth$ with
\[
\sset{f_\beta}{\beta\le\alpha}\les s_\alpha\notin \sset{s_\beta}{\beta<\alpha}.
\]
The set $\sset{s_\alpha}{\alpha<\fb}$ is $\fb$-unbounded.
\epf

By Proposition~\ref{prp:cupfin}, there are sets $X\sub\PN$ that are $\fb$-concentrated on $\Fin$.
Since $\fb$-concentrated sets are not $\sigma$-compact (Lemma~\ref{lem:nosk}),
The analogy with the resolution of Menger's Conjecture does not continue:
It will follow from our solution of the Hurewicz Problem that if $\fb=\fd$, then
there is a set $X\sub\PN$ that is $\fb$-concentrated on $\Fin$, but is \emph{not} $\ufog$.

\subsection{Omission of intervals}

We now introduce the method that we will use to settle Hurewicz's conjecture and problem.
We begin with some easy preparations.

\blem
\label{lem:x}
Let $y\in\roth$. If $k\notin y$, then $n<y(n)$ for all $n\ge k$.\qed
\elem

\bdfn
We define an operator on functions $y\in\roth$, as follows:
If $y=\bbN$ (the identity function), then $\tilde y := y$.
For $y\neq \bbN$, we define
\begin{align*}
\tilde y(1) &:= y(k), \text{ where }k:=\min y\comp\\
\tilde y(n) &:=y(\tilde y(n-1))\quad\text{for }n>1.
\end{align*}
In other words, let $k$ be the first natural number that is not in $y$.
Then $\tilde y(n)=y^n(k)$, the $n$-th iterate of the function $y$ on the number $k$.
\edfn

\begin{figure}[!htp]
	\begin{center}
		{\large
			$\xymatrix@M=0pt{
				\underset{\tilde y(1)}{\bullet}
				\ar@/^0.3cm/[rr]^y
				& &
				\underset{\tilde y(2)}{\bullet}
				\ar@/^0.3cm/[rr]^y
				& &
				\underset{\tilde y(3)}{\bullet}
				\ar@/^0.3cm/[rr]^y
				& &
				\underset{\tilde y(4)}{\bullet}
				\ar@/^0.3cm/[rr]^y
				& & &
				\dotsm
			}$
		}
	\end{center}
	\caption{The set $\tilde{y}$, for $y\neq\bbN$}
\end{figure}
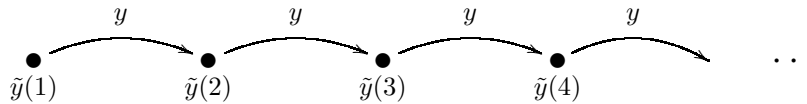

\blem
\label{lem:tilde}
Let $y\in\roth$.
\be
\item The function $\tilde y$ is increasing.
\item $y\le\tilde y$.
\ee
\elem
\bpf
This is obvious if $y=\bbN$. Assume that $y\neq\bbN$.

(1) Let $k:=\min y\comp$.
Since the function $y$ is increasing, we can iterate it
on both sides of the inequality $k<y(k)$, and have
\[
\tilde y(n)=y^{n}(k)<y^{n+1}(k)=\tilde y(n+1)
\]
for all $n$.

(2) As sets, we have $\tilde y\sub y$.
\epf

\bdfn
For natural numbers $k<m$, define the \emph{natural intervals} in the (natural) way:
\begin{align*}
\cbbl k,m \cbbr  &:= \{k,k+1,\dotsc,m\}\\
\obbl k,m \obbr   &:= \{k+1,k+2\dotsc,m-1\}\\
\cbbl k,m \obbr   &:= \{k,k+1,\dotsc,m-1\}\\
\obbl k,m \cbbr   &:= \{k+1,k+2\dotsc,m\}
\end{align*}
A set $x\in\roth$ \emph{omits} an interval $I$ if $x\cap I=\emptyset$.
\edfn

\blem
\label{lem:omit0}
Let $x,y\in\roth$, with $x\neq\bbN$.
For almost all $n$, if the set $x$ omits the interval $\ooi{\tilde y}{n}$,
then
\[
y(\tilde y(n))\le x(\tilde y(n)).
\]
\elem
\bpf
Let $k:=\min x\comp$.
For all $n\ge k$, we have $k\le y(n)$, and thus
$\tilde y(n)<x(\tilde y(n))$ (Lemma~\ref{lem:x}).
As $x(\tilde y(n))\notin \ooi{\tilde y}{n}$,
necessarily
\[
y(\tilde y(n))=\tilde y(n+1)\le x(\tilde y(n)).\qedhere
\]
\epf

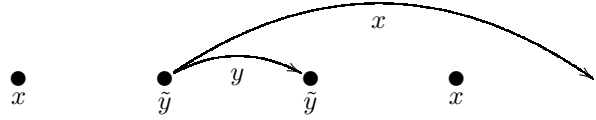
\begin{figure}[!htp]
\begin{center}
{\large
$\xymatrix@M=0pt{
\underset{x}{\bullet} & &
\underset{\tilde y}{\bullet}
\ar@/^0.3cm/[rr]_y
\ar@/^1cm/[rrrrrr]_x
& &
\underset{\tilde y}{\bullet} & &
\underset{x}{\bullet}  & &
}$
}
\end{center}
\caption{An illustration of Lemma~\ref{lem:omit0}}
\end{figure}

\bdfn
Let $a\in\roth$.
An \emph{open $a$-interval} is an interval $\ooi{a}{n}$ for some natural
number $n$. Similarly, we define \emph{closed $a$-intervals}.
\edfn

\bcor
\label{cor:omit1}
Let $x,y\in\roth$, with $x\neq\bbN$.
If the set $x$ omits infinitely many open $\tilde y$-intervals, then
$y\lei x$. \qed
\ecor

\bcor
For each function $g\in\NN$ there is a set $a\in\roth$ such that $a\comp\in\roth$ and $g\lei a,a\comp$.
\ecor
\bpf
We may assume that the function $g$ is increasing.
Take
\[
a:=\UnN{n}\cbbl\tilde g(2n-1),\tilde g(2n)\obbr,
\]
so that
\[
a\comp=\cbbl 1,\tilde g(1)\obbr\cup\UnN{n}\coi{\tilde g}{2n}.
\]
Then both sets $a$ and $a\comp$ are infinite, and they omit infinitely many
open $\tilde g$-intervals.
Apply Corollary~\ref{cor:omit1}.
\epf

\blem
\label{lem:cl}
For each family $Y\sub\roth$ with $\card{Y}<\fd$ there is a set $s\in\roth$ such that,
for each set $y\in Y$, the set $s$ omits infinitely many closed $y$-intervals.
\elem
\bpf
We may assume that the family $Y$ is closed under removing finitely many points from
its elements.
Towards a contradiction, assume that for each set $s\in\roth$ there is a set $y\in Y$ and a natural number $k$ with
\[
\card{s\cap \cci{y}{n}}\neq\emptyset
\]
for all $n\ge k$.
Removing the first $k$ elements from $y$, we may assume that this holds for all $n$.
Since the odd-numbered closed $y$-intervals are disjoint,
we have $s(n)\le y(2n)$ for all $n$.
Thus, the family of functions
$\sset{y(2n)}{y\in Y}$ is dominating.
A contradiction.
\epf

\blem
\label{lem:ld}
Let $Y\sub\roth$ be a family with $\card{Y}<\fd$.
For each set $a\in\roth$, there is a set of indices $b\in\roth$
such that the set $c:=\Un_{n\in b}\coi{a}{n}$ has $Y\lei c$.
\elem
\bpf
Let $a\in\roth$.
For a set $x\in\roth$, define
\[
x/a:=\set{n}{x\cap\coi{a}{n}\neq \emptyset}.
\]
Let $\tilde Y/a := \set{\tilde y/a}{y\in Y}$.
Take a set $b\in\roth$ such that, for each set $y\in Y$,
the set $b$ omits infinitely many closed $\tilde y/a$-intervals (Lemma~\ref{lem:cl}).
In particular, we have $b\neq\bbN$.
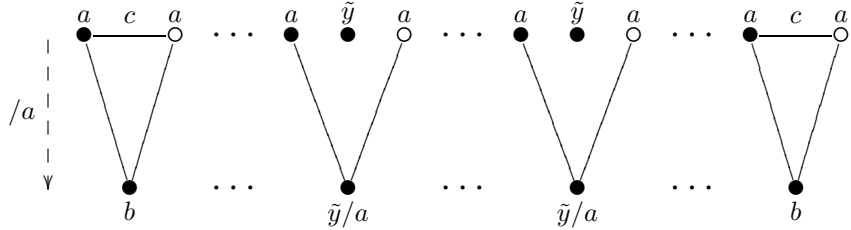
\[{\large
\xymatrix@M=0pt@R=40pt@C=10pt{
\ar@{-->}[d]_{/a}
& \overset{a}\bullet\ar@{-}[rr]^c \ar@{-}[dr]& & \overset{a}\circ \ar@{-}[dl]
& \dotsm
& \overset{a}\bullet\ar@{-}[dr] & \overset{\tilde y}\bullet &  \overset{a}\circ \ar@{-}[dl]
& \dotsm
& \overset{a}\bullet\ar@{-}[dr] & \overset{\tilde y}\bullet & \overset{a}\circ \ar@{-}[dl]
& \dotsm
& \overset{a}\bullet\ar@{-}[rr]^c\ar@{-}[dr] & & \overset{a}\circ\ar@{-}[dl]\\
& & \underset{b}\bullet & & \dotsm & &
\underset{\tilde y/a}\bullet & & \dotsm & &
\underset{\tilde y/a}\bullet & & \dotsm & & \underset{b}\bullet
}
}\]
Let
\[
c:=\Un_{n\in b}\coi{a}{n}.
\]
For each $y\in Y$, the set $c$ omits infinitely many closed $\tilde y$-intervals.
By Lemma~\ref{lem:omit0}, we have $y\lei c$.
\epf

\bcor
\label{cor:dstep}
For each family $Y\sub\roth$ with $\card{Y}<\fd$
and each set $a\in\roth$,
there is a set $s\in\roth$ with $s\comp\in\roth$, such that $Y\lei s$ and $a\lei s\comp$.
\ecor
\bpf
By Lemma~\ref{lem:ld}, there is a set $b\in\roth$ such that the set $s:=\Un_{n\in b}\coi{\tilde a}{n}$ has $Y\lei s$.
If needed, thin out the set $b$ such that its complement is also infinite.
Since the set $s\comp$ omits infinitely many
open $\tilde a$-intervals and $s\comp\neq\bbN$, we have $a\lei s\comp$.
\epf

\blem
The complement function
\begin{align*}
	\tau\colon \PN &\to\PN\\
	a &\mapsto a\comp
\end{align*}
is an isometry and, in particular, a homeomorphism.\qed
\elem

We arrive at the following solution of the Hurewicz Problem~\cite[Theorem~3.9]{sfh}.

\bthm[Tsaban--Zdomskyy]
\label{thm:HP}
There is a nontrivial $\sfoo$ set of real numbers that is not $\ufog$.
\ethm
\bpf
Fix a dominating set $\sset{f_\alpha}{\alpha<\fd}\sub\NN$.
For each ordinal $\alpha<\fd$, take
a set $s_\alpha\in\roth$ with $s_\alpha\comp\in\roth$ such that
$\sset{f_\beta}{\beta\le\alpha}\lei s_\alpha$ and $f_\alpha\lei s_\alpha\comp$
(Lemma~\ref{cor:dstep}).

The set $\sset{s_\alpha}{\alpha<\fd}$ is $\fd$-unbounded, and thus $\fd$-concentrated
on $\Fin$ (Proposition~\ref{prp:cupfin}).
By Lemma~\ref{lem:dconc}, the set $X:=\sset{s_\alpha}{\alpha<\fd}\cup\Fin$ is $\sfoo$.
The set $\tau(X)$, a continuous image of $X$ in the space $\roth$, is unbounded.
By Rec\l{}aw's Theorem~\ref{thm:rec2}, the set $X$ is not $\ufog$.
\epf

\subsection{Comments for Section~\ref{sec:hp}}
\label{com:hp}
At the correction stage of his paper~\cite[page~196]{Hure27}, Hurewicz added a proof of Sierpi\'nski
that Luzin sets are not $\ufog$. Thus, the Continuum Hypothesis settles his problem in the negative.

Within ZFC, the Hurewicz Problem was first solved by Chaber and Pol~\cite{ChaPol}, using topological methods developed by Michael.
Their proof is dichotomic:
It provides a nontrivial example if $\fb=\fd$, and a trivial one if $\fb<\fd$.
Theorem~\ref{thm:HP} provides a uniform, nontrivial example.
The author and Zdomskyy used a topological argument for Lemma~\ref{lem:ld}.
Later~\cite[Lemma~5.6]{MHP}, the author provided a combinatorial argument for this lemma,
based on one by Mildenberger~\cite{Mild01}.
Mildenberger's argument uses an ultrafilter;
the elementary proof provided here is more straightforward.
Lemma~\ref{lem:ld} is also essentially equivalent to the main result in the author's work with Zdomskyy~\cite[Theorem~4]{MGD}, whose proof was more elaborate.

The Chaber--Pol argument establishes the existence of a non-$\ufog$
set of real numbers $X$ such that all finite powers of $X$ are $\sfoo$.
The question whether there is a \emph{nontrivial} example  (of cardinality at least $\fd$)
of this kind is still open.
Chaber and Pol provide a nontrivial example when $\fb=\fd$.
The cardinal number $\fb$ is regular, that is, it is not the supremum of fewer than $\fb$ ordinals that are smaller than $\fb$.
Zdomskyy and the author proved that the regularity of the cardinal
$\fd$ suffices for a positive answer to this question~\cite[Theorem~4.6]{sfh}.

\section{The Hurewicz Conjecture}
\label{sec:hc}

\subsection{Unbounded sets}

\blem
Let $S\sub\roth$.
The set $S$ is bounded if and only if the set $\tilde S:=\sset{\tilde s}{s\in S}$ is bounded.
\elem
\bpf
$(\Leftarrow)$ By Lemma~\ref{lem:tilde}.

$(\Impl)$ Assume that $S\les b$, for some function $b\in\roth$.
For each set $s\in S$, there is a natural number $k$ such that $s\le b_k:=b\sm\cbbl 1,k \cbbr$ and, thus, $\tilde s\le \tilde b_k$.
Pick a function $c\in\roth$ with $\tilde b_k\les c$ for all $k$.
Then $\tilde S\les c$.
\epf

\blem
\label{lem:bddslalom}
Let $S\sub\roth$. The following assertions are equivalent:
\be
\itm The set $S$ is unbounded.
\itm For each set $a\in\roth$, there is a set $s\in S$ that omits infinitely many open $a$-intervals.
\ee
\elem
\bpf
$(2)\Impl(1)$: Let $a\in\roth$.
By thinning out $a$, we may assume that all open
$a$-intervals, and thus all open $\tilde a$-intervals, are nonempty.
Take a set $s\in S$ that omits infinitely many open $\tilde a$-intervals.
In particular, we have $s\neq\bbN$.
Then $a\lei s$ (Corollary~\ref{cor:omit1}).

$(1)\Impl(2)$: Assume that there is a set $a\in\roth$ such that for each set $s\in S$, we have
\[
s\cap\ooi{a}{n}\neq\emptyset
\]
for almost all $n$.
Let $k$ be a natural number with $s\cap\ooi{a}{n}\neq\emptyset$ for all $n\ge k$.
Then $s(n) < a_k(n):=a(n+k)$ for all $n$.
Let $b\in\NN$ be a function with $a_k\les b$ for all $k$. Then $S\les b$.
\epf

\bdfn
A \emph{$\fb$-scale} is an unbounded set $\sset{s_\alpha}{\alpha<\fb}\sub\roth$ such that
the enumeration is increasing with respect to $\les$, that is, we have
$s_\alpha\les s_\beta$ for all ordinals $\alpha<\beta<\fb$.
\edfn

\blem
There are $\fb$-scales.
\elem
\bpf
Let $\sset{f_\alpha}{\alpha<\fb}\sub\NN$ be unbounded.
For each ordinal $\alpha<\fb$,
take a function $s_\alpha\in\roth$ with
\[
\sset{s_\beta}{\beta<\alpha}\cup\{f_\alpha\}\les s_\alpha.\qedhere
\]
\epf

\blem
\label{lem:bsc-bu}
Every $\fb$-scale is $\fb$-unbounded.
\elem
\bpf
Let $S=\sset{s_\alpha}{\alpha<\fb}\sub\roth$ be a $\fb$-scale.
Let $g\in\NN$.
Since the set $S$ is unbounded, there is an ordinal $\alpha$ with $g<^\oo s_\alpha$.
For each ordinal $\beta\ge \alpha$ we have
\[
g <^\oo s_\alpha \les s_\beta,
\]
and thus $g <^\oo s_\beta$.
\epf

\subsection{Omission of intervals and open sets}

\bdfn
An open cover $\cU$ of a topological space $X$ is an
\emph{$\w$-cover} if every finite set $F\sub X$ is contained in some member of the cover.
Let $\Om(X)$ be the family of all $\w$-covers of the space $X$.
\edfn

Every point-cofinite cover of a space is an $\w$-cover of that space.
The converse is not true in general, but being an $\w$-cover may help establishing being a point-cofinite cover.

\blem\mbox{}
\label{lem:om}
\be
\item If $\cU$ is an $\w$-cover of a space $X$, then for each finite set $F\sub X$,
the set $F$ is contained in infinitely many members of the cover.
\item If a family $\cU=\sset{U_n}{n\in\bbN}$ (where the sets $\seq{U}$ are not necessarily distinct)
has no element covering the entire set $X$, and every point $x\in X$ is in $U_n$ for almost all $n$, then
the family $\cU$ is a point-cofinite cover of the set $X$.
%
%
\ee
\elem
\bpf
(1) Recall that we assume that $X\notin\cU$ (Convention~\ref{cnv:cover}).

Pick a set $U_1\in\cU$ with $F\sub U_1$, and a point $x_1\notin U_1$.
For each natural numbers $n>1$,
pick a set $U_n\in\cU$ with
\[
F\cup\{x_1,\dotsc,x_{n-1}\}\sub U_n
\]
and a point $x_n\notin U_n$.

(2) The family $\cU=\sset{U_n}{n\in\bbN}$ is an $\w$-cover of $X$. By (2), it is infinite.
\epf

The following observation~\cite[Proof of Lemma~1.2]{GM84} is one of the main tools of omission of intervals.

\blem[Galvin--Miller]
\label{lem:gm}
Assume that $\Fin\sub X\sub \PN$.
For each cover $\cU\in\Om(X)$,
there are a function $a\in\roth$ with $a(1)=1$ and distinct sets $\seq{U}\in\cU$
such that, for each $n$,
all points $x\in\PN$ that omit the interval $\ooi{a}{n}$
are in the set $U_n$.
\elem
\bpf
Let $a(1):=1$.
Take an open set $U_1\in\cU$ with $\emptyset,\{1\}\in U_1$.
Since the set $U_1$ is open, there is a natural number $a(2)>1$ such that
every set $x\in\PN$ with
\[
x\cap \cbbl 1,a(2) \obbr\in\{\, \emptyset,\{1\}\, \}
\]
is in the set $U_1$.
If a set $x\in\PN$ omits the interval $\obbl a(1),a(2)\obbr$,
then $x\cap \cbbl 1,a(2)\obbr\in\{\, \emptyset,\{1\}\, \}$, and thus
$x\in U_1$.

For each natural number $n>1$:
Since $\cU\in\Om(X)$, the finite set $\op{P}( \cbbl 1,a(n)\cbbr )$ is contained in
infinitely many members of the cover $\cU$.
Take an open set $U_n\in\cU\sm\{U_1,\dotsc,U_{n-1}\}$
such that $\op{P}( \cbbl 1,a(n)\cbbr )\sub U_n$.
As the set $U_n$ is open, there is a natural number $a(n+1)>a(n)$ such that
every set $x\in\PN$ with
\[
x\cap \cbbl 1,a(n+1)\obbr\in \op{P}(\cbbl 1,a(n)\cbbr)
\]
is in the set $U_n$.

Let $n$ be a natural number.
If a set $x\in\PN$ omits the interval $\ooi{a}{n}$,
then
\[
x\cap \cbbl 1,a(n+1)\obbr\sub \cbbl 1,a(n)\cbbr,
\]
and thus $x\in U_n$.
\epf

\nc\gmline[2]
{
				\underset{#1(1)}{\bullet}
				\ar@/^0.3cm/[rr]^{\gm{#2_1}}
				& &
				\underset{#1(2)}{\bullet}
				\ar@/^0.3cm/[rr]^{\gm{#2_2}}
				& &
				\underset{#1(3)}{\bullet}
				\ar@/^0.3cm/[rr]^{\gm{#2_3}}
				& &
				\underset{#1(4)}{\bullet}
				\ar@/^0.3cm/[rr]^{\gm{#2_4}}
				& &
				\underset{#1(5)}{\bullet}
				\ar@/^0.3cm/[rr]^{\gm{#2_5}}
				& & &
				\dotsm
}

\begin{figure}[!htp]
\begin{center}
{\large
$\xymatrix@M=0pt{
\gmline{a}{U}
}$
}
	\end{center}
	\caption{An illustration of the Galvin--Miller Lemma}
\end{figure}
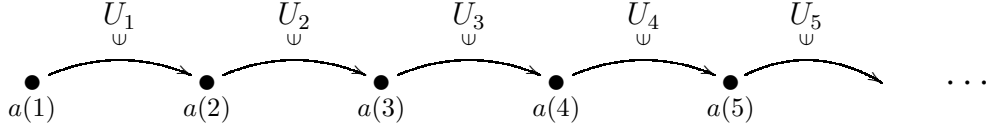

The proof of the following lemma~\cite[Theorem~11]{coc2} is similar to that of Lemma~\ref{lem:sfgo}.

\blem[Scheepers]
\label{lem:ufgg}
For Lindel\"of spaces, $\ufog=\ufin(\Ga,\Ga)$. \qed
\elem

By Lemma~\ref{lem:bsc-bu}, Proposition~\ref{prp:cupfin} and Lemma~\ref{lem:nosk}, the following theorem~\cite[Theorem~1]{BaSh01}
settles the Hurewicz Conjecture.

\bthm[Bartoszy\'nski--Shelah]
\label{thm:bs}
For each $\fb$-scale $S$, the set $S\cup\Fin$ is $\ufog$.
\ethm
\bpf
Let $S := \sset{s_\alpha}{\alpha<\fb}\sub\roth$ be a $\fb$-scale.
It suffices to prove that the set $S\cup\Fin$ is $\ufin(\Ga,\Ga)$.

Let $\seq\cU\in\Ga(X)$.
For each $n$, take a function $a_n\in\roth$ and distinct sets $\seq{U^n}$ for the cover $\cU_n$ as in Lemma~\ref{lem:gm}.
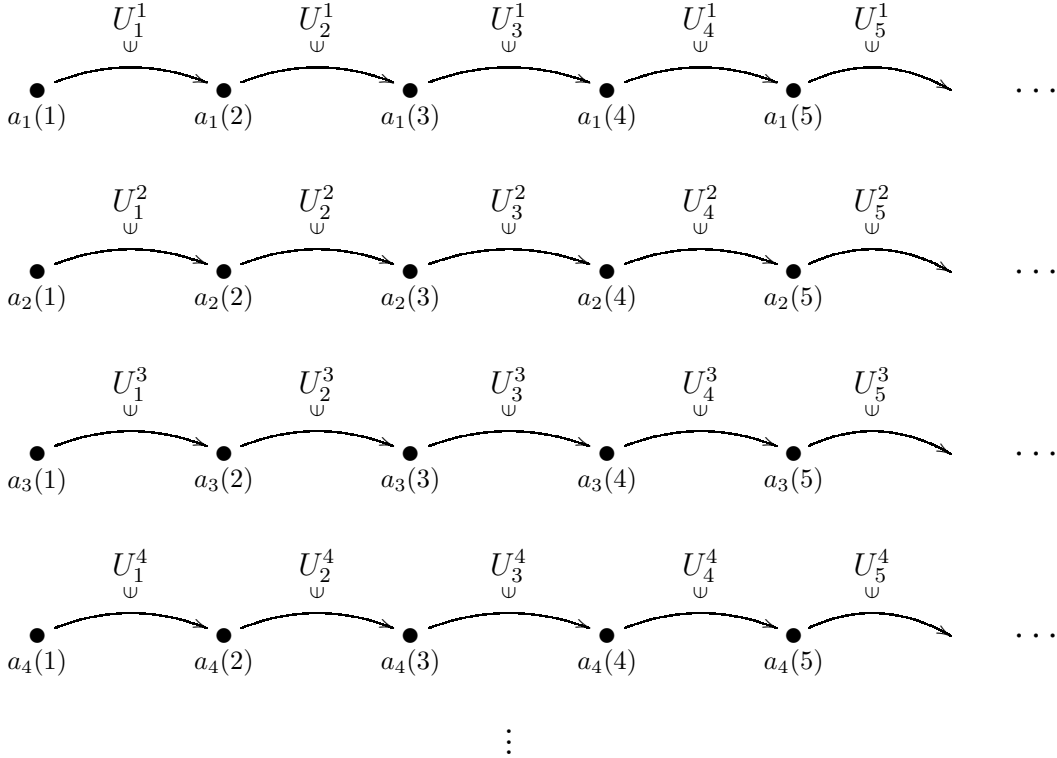
\begin{center}
{\large
$\xymatrix@M=0pt{
\gmline{a_1}{U^1}\\
\\
\gmline{a_2}{U^2}\\
\\
\gmline{a_3}{U^3}\\
\\
\gmline{a_4}{U^4}\\
& & & & & \vdots
}$
}
\end{center}
Pick an ordinal $\alpha<\fb$ such that the set
\[
I:=\sset{n}{a_n(n+1)<s_\alpha(n)}
\]
is infinite.
By Corollary~\ref{cor:nonh}, the set $A := \Fin\cup\sset{s_\beta}{\beta<\alpha}$ is $\ufin(\Op,\Ga)$.
By adding countably many points to the set $A$, we may assume that no given cover contains a
finite cover of the set $A$ (we add witnesses that there are no finite subcovers of the set $X$).
Thus, there are finite sets $\cG_n\sub\cU_n$, for $n\in I$, such that
\[
\sset{\Un \cG_n}{n\in I}\in\Ga(\Fin\cup\sset{s_\beta}{\beta<\alpha}).
\]
Let $\cF_n:=\emptyset$ for $n\nin I$, and
\[
\cF_n:=\{U^n_1,\dotsc,U^n_n\}\cup\cG_n
\]
for $n\in I$.
Then
\[
\sset{\Un\cF_n}{n\in\bbN}=\sset{\Un\cF_n}{n\in I}\cup\{\emptyset\}.
\]
By Lemma~\ref{lem:om}(2), it suffices to show that for each point $x\in X$, we have $x\in\Un\cF_n$ for almost all $n\in I$.

If $x\in \Fin\cup\sset{s_\beta}{\beta<\alpha}$, this is evident.
Fix an ordinal $\beta\ge\alpha$.
Then $s_\alpha \les s_\beta$.
It follows that for almost all $n\in I$:
\[
a_n(n+1)< s_\alpha(n) \le s_\beta(n),
\]
and thus $\card{s_\beta\cap \cbbl 1,a_n(n+1)\obbr}<n$.
The set $\cbbl 1,a_n(n+1)\obbr$ contains exactly $n$ disjoint
open $a_n$-intervals.
Thus, there is $i\le n$ such
the set $s_\beta$ omits the interval $\ooi{a_n}{i}$,
and thus $s_\beta\in U^n_i\sub\Un\cF_n$.
\epf

\subsection*{Comments for Section~\ref{sec:hc}}
The Continuum Hypothesis implies the existence of a \emph{Sierpi\'nski set}:
An uncountable set $S$ of real numbers whose intersection with every Lebesgue null set is countable.
Sierpi\'nski sets are $\sogg$~\cite[Theorems~2.9--2.10]{coc2} and, in particular, are $\ufog$ and not $\sigma$-compact (Proposition~\ref{prp:nosk}).

Just, Miller, Scheepers and Szeptycki~\cite[Theorem~5.1]{coc2} were the first to settle the Hurewicz Conjecture~\ref{cnj:Hur} in ZFC.
Here, too, the argument was dichotomic ($\aleph_1=\fb$ or $\aleph_1<\fb$), and one of the cases was trivial.
The Bartoszy\'nski--Shelah Theorem~\ref{thm:bs} provided the first uniform, nontrivial counterexample, using
a topological proof.
The proof provided here, using omission of intervals, is a modification of an earlier proof by the author~\cite[Theorem~2.12]{MHP}.
The notion of $\fb$-scale in Theorem~\ref{thm:bs} can be generalized~\cite[Theorem~5.4]{pMReal};
for example, the construction in the proof of Lemma~\ref{lem:b-unb} suffices.
We will use the flexibility of our proof below.

In the comments to Section~\ref{sec:hp}, we asked whether, provably, there is a nontrivial $\sfoo$ set of real numbers
in all finite powers that is not $\ufog$.
The analogous question for $\ufog$ was solved by Bartoszy\'nski and the author~\cite[Theorem~10]{ideals}:
For each $\fb$-scale $S$, the non-$\sigma$-compact set $S\cup\Fin$ is $\ufog$ in all finite powers.
Zdomskyy and the author proved that, in fact,
the product $(S\cup\Fin)\x Y$ is $\ufog$ for every $\ufog$ set of real numbers $Y$~\cite[Theorem~2.28]{sfh}.
Szewczak and the author~\cite[Theorem~5.4]{pMReal} generalized this result further.

\section{Almost uniform convergence of subsequences}
\label{sec:sogg}

\subsection{Diagonalizing point-cofinite covers}
For a topological space $X$, let $\Cp(X)$ be the family of continuous real-valued functions on the space $X$,
endowed with the topology of pointwise convergence of nets, or equivalently, viewed as a subspace of the Tychonoff product space $\bbR^X$.
For functions $f,\seq{f}\in\Cp(X)$ we have $f_n\longrightarrow f$ if and only if $f_n(x)\to f(x)$ for every point $x\in X$.

The following notions are motivated by old studies in trigonometric series;
we refer to the introduction of Bukovsk\'y, Rec\l{}aw and Repick\'y for details~\cite{BRR91}.
A convergence of real-valued functions $f_n\longrightarrow f$ is \emph{quasinormal} if there is a sequence of positive real numbers $\epsilon_n\longrightarrow 0$
such that for each point $x\in X$ we have
\[
\card{f_n(x)-f(x)}\le \epsilon_n
\]
for almost all $n$.
A topological space $X$ is  a \emph{QN space} if every convergent sequence in the space $\Cp(X)$ converges quasinormally.
It is a \emph{wQN space} if every convergent sequence in $\Cp(X)$ has a subsequence that converges quasinormally.

Scheepers~\cite[Theorem~8]{wqn} proved that every $\sone(\Ga,\Ga)$ space is wQN.
Sakai~\cite[Theorem~2.5]{Sakai07} and Bukovsk\'y--Hale\v{s}~\cite[Theorem~9]{BH07} proved that a subset of the Cantor space is wQN if and only if it is $\sone(\Gamma,\Gamma)$ for clopen covers.
Since $\ufin(\Ga,\Ga)=\ufog$, the property $\sogg$ implies $\ufog$.
By Proposition~\ref{prp:nosk}, the property $\sogg$ is strictly stronger than $\ufog$.
However, the critical cardinalities of these properties are equal.

\blem
\label{lem:nons1gg}
$\non(\sogg)=\fb$.
\elem
\bpf
$(\le)$: By Corollary~\ref{cor:nonh}.

$(\ge)$: Suppose that $\card{X}<\fb$ and $\seq{\cU}\in\Ga(X)$.
By moving to countably infinite subsets, we can enumerate
\[
\cU_n=\set{U^n_m}{m\in\bbN}
\]
for each natural number $n$.
For each point $x\in X$, define a function $f_x\in\NN$ by
\[
f_x(n):=\min\set{m}{x\in U^n_k \text{ for all }k\ge m}.
\]
Let $g\in\NN$ be a $\les$-upper bound for all these functions $f_x$.
Then $\set{U^n_{g(n)}}{n\in\bbN}\in\Ga(X)$.
\epf

\bdfn
For sets $a,b\in\roth$, we write $a\subs b$ if the set $a\sm b$ is finite, that is, if $a(n)\in b$ for almost all $n$.
\edfn

\bdfn
Let $\kappa$ be an uncountable cardinal number.
A $\kappa$-\emph{tower} is a set $\sset{t_\alpha}{\alpha<\kappa}\sub\roth$ such that
$t_\beta\subs t_\alpha$ for all $\alpha<\beta$.
A set is a \emph{tower} if it is a $\kappa$-tower for some cardinal $\kappa$.
\edfn

\bdfn
For sets $a,b\in\roth$, we write $a\sqe b$ if almost all closed $b$-intervals contain closed $a$-intervals.
Equivalently, if
\[
\card{a\cap\cci{b}{n}}\ge 2
\]
for almost all $n$.
\edfn

The proof of the following lemma is similar to that of Lemma~\ref{lem:bddslalom}.

\blem
\label{lem:bdd2}
A set $S\sub\roth$ is $\les$-bounded if and only if it is $\sqe$-bounded.\qed
\elem

Since a tower is a subset of $\roth$, it may be unbounded.
The following theorem is implicit in a proof of Scheepers~\cite[Theorem~6]{Scheepers98},
and explicit in Miller--Tsaban~\cite[Proposition~2.5]{BBC}.

\bthm
\label{thm:btowersogg}
For each unbounded $\fb$-tower $T$, the set $T\cup\Fin$ is $\sogg$.
\ethm
\bpf
Let $T=\sset{t_\alpha}{\alpha<\fb}$ be an unbounded $\fb$-tower.
Let $\seq{\cU}\in\Ga(T\cup\Fin)$.
For each $n$, take a function $a_n\in\roth$ and distinct sets $\seq{U^n}$ for the cover $\cU_n$ as in Lemma~\ref{lem:gm}.
The set $\sset{a_n}{n\in\bbN}$ is countable and thus bounded.
By Lemma~\ref{lem:bdd2}, there is a set $a\in\roth$ with $a_n\sqe a$ for all $n$.

Since the set $T$ is unbounded, there is an ordinal $\alpha<\fb$ such that
$t_\alpha$ omits infinitely many open $a$-intervals (Lemma~\ref{lem:bddslalom}).
It follows that for each $n$, the set $t_\alpha$ omits infinitely many
open $a_n$-intervals.
For each $n$, let
\[
I_n := \sset{m\ge n}{t_\alpha\cap \ooi{a_n}{m}=\emptyset}.
\]
Being an infinite subset of $\cU_n$, we have $\sset{U^n_m}{m\in I_n}\in\Ga(T\cup\Fin)$.

By Lemma~\ref{lem:nons1gg}, there is for each $n$ a natural number $g(n)\in I_n$ such that
\[
\sset{ U^n_{g(n)} }{n\in\bbN}\in\Ga(\Fin\cup\sset{t_\beta}{\beta<\alpha})
\]
For $\beta\ge\alpha$, we have $t_\beta\subs t_\alpha$.
Let $k$ be a natural number such that $t_\beta\sm\cbbl 1,k\cbbr\sub t_\alpha$.
For each $n\ge k$, we have $k\le n \le g(n)$, and thus $k\le a_n(g(n))$.
It follows that
\[
t_\beta\cap \ooi{a_n}{g(n)} \sub t_\alpha\cap \ooi{a_n}{g(n)}=\emptyset.
\]
Consequently, we have $t_\beta\in U^n_{g(n)}$.
It follows that $\sset{ U^n_{g(n)} }{n\in\bbN}\in\Ga(T\cup\Fin)$.
\epf

Theorem~\ref{thm:btowersogg} is fundamentally different from the previous theorems:
The existence of an unbounded $\fb$-tower is not provable (in ZFC).
Even the mere existence of a nontrivial $\sogg$ set of real numbers is not provable~\cite[Theorem~3.2]{BBC}.
If, e.g., $\aleph_1=\fb$ (Lemma~\ref{lem:ptower} below) or $\fb<\fd$ (see the following proof), then there is an unbounded $\fb$-tower.
We conclude this section by pointing out that if there is an unbounded tower (of any height), then there is an unbounded $\fb$-tower~\cite[Proposition~2.4]{BBC}.

\bprp[Miller--Tsaban]
\label{prp:anytower}
If there is an unbounded tower, then there is an unbounded $\fb$-tower.
\eprp
\bpf
The proof is dichotomic.

Assume that $\fb<\fd$.
Let $S=\sset{s_\alpha}{\alpha<\fb}\sub\NN$ be a $\fb$-scale.
As $\card{S}<\fd$, there is a function $g\in\NN$ with $S\lei g$.
For each ordinal $\alpha<\fb$, let
\[
t_\alpha := \sset{n}{s_\alpha(n)\le g(n)}.
\]
Clearly, $t_\beta\subs t_\alpha$ for $\alpha<\beta$.

Let $T:=\sset{t_\alpha}{\alpha<\fb}$.
Assume that $T\les h$ for some function $h\in\NN$.
For each $\alpha<\fb$, we have
\[
s_\alpha(n)\le s_\alpha(t_\alpha(n))\le g(t_\alpha(n))\le g(h(n))
\]
for almost all $n$. Thus, $T\les g\circ h$.
A contradiction.

Next, assume that $\fb=\fd$.
Let $T$ be an unbounded tower.
Fix a dominating family $\sset{f_\alpha}{\alpha<\fb}\sub\NN$.
For each ordinal $\alpha<\fb$,
let $g$ be a $\les$-bound of the set $\sset{t_\beta\sm \cbbl 1,k\cbbr}{\beta<\alpha, k\in\bbN}$,
and pick an element $t_\alpha\in T$ with $f_\alpha, g <^\oo t_\alpha$.
Then
\[
\{ f_\alpha \} \cup \sset{t_\beta\sm \cbbl 1,k\cbbr}{\beta<\alpha, k\in\bbN} <^\oo t_\alpha.
\]
The set $\sset{t_\alpha}{\alpha<\fb}$ is unbounded.
Let $\alpha<\beta<\fb$ be ordinals.
Towards a contradiction, assume that $t_\alpha\subs t_\beta$.
Then there is $k$ with $t_\alpha\sm \cbbl 1,k\cbbr\sub t_\beta$, so that
\[
t_\beta \le t_\alpha\sm \cbbl 1,k\cbbr <^\oo t_\beta.
\]
A contradiction. It follows that $t_\beta\subs t_\alpha$.
\epf

\bcor[Miller--Tsaban]
\label{cor:nontsogg}
If there is an unbounded tower (of any cardinality), then
there is a nontrivial $\sogg$ set of real numbers.
\ecor

\subsection{Taking unions of pairs of elements}

A seemingly minor relaxation of the property $\sogg$ brings us back to ZFC.

\bdfn
Let $\rmA$ and $\rmB$ be types of covers of topological spaces.
A topological space $X$ is $\mathsf{U}_2(\rmA,\rmB)$ if for each sequence of covers
$\seq{\cU}\in\rmA(X)$, none with a two-element subcover,
there are sets $U_1,V_1\in\cU_1$, $U_2,V_2\in\cU_2$, \dots
such that
\[
\{U_1\cup V_1, U_2\cup V_2,\dotsc\}\in\rmB(X).
\]
\edfn

By Lemma~\ref{lem:ufgg}, we have the following implications
\[
\sogg \longrightarrow \utgg \longrightarrow \ufog.
\]
We will later see that these implications are strict.
Now, we adjust our methods to produce $\utgg$ sets.

\bdfn
A \emph{$\fb(\sqe)$-scale} is an unbounded set $\sset{s_\alpha}{\alpha<\fb}\sub\roth$
such that $s_\alpha\sqe s_\beta$ for all ordinals $\alpha<\beta<\fb$.
\edfn

\blem
\label{lem:bscip}
There are $\fb(\sqe)$-scales.
\elem
\bpf
Let $\sset{f_\alpha}{\alpha<\fb}\sub\NN$ be an unbounded set.
By Lemma~\ref{lem:bdd2},
for each ordinal $\alpha<\fb$,  there is
a set $s\in\roth$ with $\sset{s_\beta}{\beta<\alpha}\sqe s$.
Take an infinite set $s_\alpha\sub s$ with $f_\alpha\les s_\alpha$.
Then $\sset{s_\beta}{\beta<\alpha}\sqe s_\alpha$.
\epf

\blem
\label{lem:gm2}
Assume that $\Fin\sub X\sub \PN$.
For each cover $\cU\in\Om(X)$,
there are a function $b\in\roth$ with $b(1)=1$ and distinct sets $\seq{U}\in\cU$
such that, for each $n$,
all points $x\in\PN$ with
\[
\card{x\cap\ooi{b}{n}}\le 1,
\]
are in the set $U_{2n-1}\cup U_{2n}$.
\elem
\bpf
Take a set $a\in\roth$ and distinct sets
$\seq{U}$ for the cover $\cU$ as in Lemma~\ref{lem:gm}.
Let $b:=\sset{a(2n-1)}{n\in\bbN}$.
If a set $x\in\PN$ has
\[
\card{x\cap\ooi{b}{n}}\le 1,
\]
then the set $x$ omits at least one of the two $a$-intervals,
$\ooix{a}{2n-1}{2n}$ and $\ooix{a}{2n}{2n+1}$
and thus $x\in U_{2n-1}\cup U_{2n}$.
\epf

We now see that there are nontrivial $\utgg$ sets of real numbers in ZFC~\cite[Theorem~4.4]{MHP}.

\bthm[Tsaban]
\label{thm:utgg}
For each $\fb(\sqe)$-scale $S$, the set $S\cup\Fin$ is $\utgg$.
\ethm
\bpf
Let $S=\sset{s_\alpha}{\alpha<\fb}$ be a $\fb(\sqe)$-scale.
Let $\seq{\cU}\in\Ga(S\cup\Fin)$.

By Lemma~\ref{lem:gm2},
for each $n$ there are a function $b_n\in\roth$ with $b_n(1)=1$ and distinct sets $\seq{U^n}\in\cU_n$
such that, for each $m$,
all points $x\in\PN$ with
\[
\card{x\cap\ooi{b_n}{m}}\le 1,
\]
are in the set
\[
U^n_{2m-1}\cup U^n_{2m}.
\]
By Lemma~\ref{lem:bdd2}, there is a set $b\in\roth$ with $\sset{b_n}{n\in\bbN}\sqe b$.
There is a function $c\in\roth$ such that, for each $n$, for all $m\ge c(n)$,
the interval $\ooi{b}{m}$ contains an open $b_n$-interval,
say $\ooi{b_n}{i_{n,m}}$.
If $m\ge c(n)$ and a set $x\in\PN$ has
\[
\card{x\cap\ooi{b}{m}}\le 1,
\]
then $\card{x\cap\ooi{b_n}{i_{n,m}}}\le 1$,
and thus
\[
x\in W^n_m := U^n_{2i_{n,m}-1}\cup U^n_{2i_{n,m}}.
\]
As the set $S$ is unbounded, there is by Lemma~\ref{lem:bddslalom} an ordinal $\alpha<\fb$,
such that the set
\[
I := \sset{m}{s_\alpha \cap \ooi{b}{m} = \emptyset}
\]
is infinite.
Since the set $I$ is infinite, for each $n$ we have
\[
\set{W^n_{m}}{c(n)\le m\in I}\in\Ga(S\cup\Fin).
\]
The set $\Fin\cup\set{s_\beta}{\beta<\alpha}$ is $\sogg$ (Lemma~\ref{lem:nons1gg}).
Thus, there is a function $g\in\NN$ with $c(n)\le g(n)\in I$ for all $n$, such that
\[
\set{W^n_{g(n)}}{n\in\bbN}\in\Ga(\Fin\cup\set{s_\beta}{\beta<\alpha}).
\]
We have $n\le c(n)\le g(n)$ for all $n$.
Let $\beta\ge \alpha$. Then $s_\alpha\sqe s_\beta$.
Let $k$ be a natural number such that every open interval $\ooi{s_\beta}{m}$ with $s_\beta(m)\ge k$ contains an open $s_\alpha$-interval.
For $n\ge k$, since $g(n)\in I$,
the set $s_\alpha$ omits the interval $\ooi{b}{g(n)}$.
Suppose that the set $s_\beta$ has more than one point in the interval $\ooi{b}{g(n)}$,
that is, the interval $\ooi{b}{g(n)}$ contains a \emph{closed} $s_\beta$-interval.
Since $k\le n\le g(n)\le b(g(n))$, this closed $s_\beta$-interval contains a closed $s_\alpha$-interval.
A contradiction.
Thus, for almost all $n$, we have $\card{s_\beta\cap\ooi{b}{g(n)}}\le 1$ and thus $s_\beta\in W^n_{g(n)}$.

For each $n$, the set $W^n_{g(n)}$ is a union of two sets from the cover $\cU_n$.
\epf

\subsection{Intermediate properties}
\label{subsec:inter}

\bdfn
Let $\rmA$ and $\rmB$ be types of covers of topological spaces, and let $f\in\NN$.
A topological space $X$ is $\mathsf{U}_f(\rmA,\rmB)$ if
for each sequence
$\seq{\cU}\in\rmA(X)$, none with a finite subcover,
there are subsets $\cF_1\sub\cU_1$, $\cF_2\sub\cU_2$, \dots, with $\card{\cF_n}\le f(n)$ for all $n$,
such that $\{\seq{\Un\cF}\}\in\rmB(X)$.

If $f$ is a constant function $f(n)=k$ for all $n$, we write $\mathsf{U}_k(\rmA,\rmB)$ instead of
$\mathsf{U}_f(\rmA,\rmB)$.
\edfn

We have the following observation~\cite[Lemma~3.2]{MHP}.
\blem
\label{lem:limsup}
For families $\rmB$ such that $\{\emptyset\}\cup\cU\in\rmB$ for all $\cU\in\rmB$,
the property $\mathsf{U}_f(\rmA,\rmB)$ depends only on $\limsup_n f(n)$.
\elem

We thus have the following infinite hierarchy, where $\op{id}\in\NN$ is the identity function:
\[
\sogg \longrightarrow \utgg \longrightarrow \mathsf{U}_3(\Ga,\Ga) \longrightarrow \dotsm \longrightarrow \uidgg \longrightarrow \ufog.
\]
The property $\uidgg$ is strictly stronger than $\ufog$.
For example, it implies $\sone(\Ga,\Op)$~\cite[Theorem~3.5]{MHP}, a property that Cantor's space does not have (Proposition~\ref{prp:nosk}).
We have already mentioned that the properties $\sogg$ and $\utgg$ are not equivalent.
We will later see that all of these implications are strict.
The following result~\cite[Theorem~3.3]{MHP} improves upon Theorem~\ref{thm:bs}.

\bthm[Tsaban]
\label{thm:un}
For each $\fb$-scale $S$, the set $S\cup\Fin$ is $\uidgg$.
\ethm
\bpf
In our proof of Theorem~\ref{thm:bs}, we apply $\sogg$ instead of $\ufog$ to cover the set of cardinality smaller than $\fb$.
Thus, it suffices to choose $n+1$ sets.
\epf

\subsection{Comments for Section~\ref{sec:sogg}}

The \emph{Scheepers Conjecture}~\cite[Conjecture~1]{wqn} asserts that
for subsets of the Cantor space, the property $\sogg$ restricted to clopen covers is equivalent to $\sogg$ (for open covers).
After two decades of study of this conjecture by various mathematicians, Peng~\cite[Section~4]{Peng} proved that this conjecture is false:
Assuming \CH{}, there is a subset of the Cantor space that is $\sogg$ for clopen covers but not for arbitrary open covers.
If we consider countable Borel covers instead of open covers, the properties corresponding to $\sogg$, wQN and $\ufog$ are all equivalent,
and constitute a much stronger property~\cite{cbc, hH}.

The property wQN is not the only manifestation of $\sogg$ (for clopen covers) in the realm of function spaces.
Arhangel'ski\u{\i}~\cite{Arhan72, Arhan81} introduced four properties of topological space, that
help determining when a product of Fr\'echet--Urysohn spaces is Fr\'echet--Urysohn.
He named these properties $\alpha_1$ to $\alpha_4$.
A topological space is $\alpha_2$ if for every sequence of sequences converging to the same element,
we can choose one element from each sequence, and obtain a sequence that converges to the same element.
Scheepers~\cite[Theorems~8 and~9]{wqn}, based on Fremlin~\cite[Section~12]{FremlinSeq}, proved that for a topological space $X$,
the function space $\Cp(X)$ is $\alpha_2$ if and only if the space $X$ is wQN.

The proof of Proposition~\ref{prp:anytower} bridges a minor gap in the author's original argument.

Miller and the author~\cite[Theorem~2.8]{BBC} proved that if the union of fewer than $\fb$ set that are $\sogg$ is $\sogg$ then for each unbounded $\fb$-tower $T$, all finite powers of the set $T\cup\Fin$ are $\sogg$.
Szewczak and W\l{}udecka~\cite[Theorem~3.1]{SzWl21} prove the same assertion \emph{without any hypothesis}.
Moreover, they establish that every product of such sets is $\sogg$, and even of more general types of sets.

The original Bartoszy\'nski--Shelah proof of Theorem~\ref{thm:bs} does not imply Theorem~\ref{thm:un}.
Theorem~\ref{thm:utgg} was first established for a slightly more general family than $\fb(\sqe)$-scales~\cite[Theorem~4.4]{MHP},
but the relation used here is transitive, which is convenient.

\section{The Crown Theorem}
\label{sec:crown}

Recall that for a topological space $X$, the space $\Cp(X)$ consists of the continuous real-valued functions on $X$, with the topology of pointwise convergence.
The basic open neighborhoods of the zero function are the sets
\[
[F,m] := \smallmedset{f\in\Cp(X)}{ \card{f(x)} < \tfrac{1}{m} \text{ for } x\in F },
\]
for $F\sub X$ finite and $m\in\bbN$.

A topological space $X$ is Tychonoff if for each point $a$ and every closed set $C$ with $a\notin C$,
there is a function $f\in\Cp(X)$ that is $0$ on the set $C$ and $1$ on the point $a$.
Every metric space, in particular every set of real numbers, is Tychonoff.
We have the following fact.

\bprp
For an infinite Tychonoff space $X$, the space $\Cp(X)$ is metrizable if and only if the set $X$ is countable.
\eprp
\bpf
$(\Leftarrow)$ If the set $X$ is countable, then the space $\bbR^X$ is homeomorphic to the metrizable space $\bbR^\bbN$.
The space $\Cp(X)$ is a subspace of $\bbR^X$.

$(\Impl)$ Suppose that the space $\Cp(X)$ is metrizable.
In particular, there is a countable local base at the zero function $0$.
We may assume that it is of the form
$\sset{[F_n,m_n]}{n\in\bbN}$.
Let $a\in X\sm \UnN{n}{F_n}$.
For each $n$, we have $[F_n,m_n]\nsubseteq [\{a\},1]$:
Let $U$ be an open neighborhood of the point $a$ that is disjoint from the finite set $F_n$.
Since the space $X$ is Tychonoff, there is a function $f\in\Cp(X)$ that is $0$ on the set $U\comp$ and $1$ on the point $a$.
\epf

A topological space $Y$ is \emph{Fr\'echet--Urysohn} if for each set $A\sub Y$,
every point in the closure of $A$ is the limit of a convergent sequence of points from the set $A$.
Every metric space is Fr\'echet--Urysohn, and it is natural to ask whether, for uncountable sets of real numbers $X$,
the function space $\Cp(X)$ may be Fr\'echet--Urysohn.

\bdfn
Let $\rmA$ and $\rmB$ be types of covers of topological spaces.
A topological space $X$ is $\smallbinom{\rmA}{\rmB}$ if every cover in $\rmA(X)$ has a subcover in $\rmB(X)$.
\edfn

For example, \emph{compact} and \emph{Lindel\"of} are properties of this kind.
The following result~\cite[Theorem~2]{GN82} translates the Fr\'echet--Urysohn property of the space $\Cp(X)$ to a covering property of $X$.

\bthm[{Gerlits--Nagy}]
Let $X$ be a set of real numbers.
The space $\Cp(X)$ is Fr\'echet--Urysohn if and only if the set $X$ is $\omg$.
\ethm

To carry out their proof, Gerlits and Nagy establish the following result.

\blem
$\omg=\sone(\Om,\Ga)$.
\elem
\bpf
$(\Leftarrow)$ Let $\cU\in\Om(X)$.
Apply $\sone(\Om,\Ga)$ to the constant sequence
$\cU,\cU,\cU,\dotsc$ to obtain elements
$\seq{U}\in\cU$ with $\{\seq{U}\}\in\Ga(X)$.

$(\Impl)$ Assume that $X$ is $\omg$.
Let $\seq{\cU}\in\Om(X)$.
For each $n$, we have
\[
\cU_n' := \set{U_1\cap\dotsb\cap U_n}{U_1\in\cU_1,\dotsc,U_n\in\cU_n}\in\Om(X),
\]
the cover $\cU_n'$ refines the cover $\cU_n$, and the covers $\cU_n'$ get finer as $n$ increases.
By Lemma~\ref{lem:om}(2), if we enlarge the open sets in a point-cofinite cover the cover remains point-cofinite.
Thus, it suffices to consider these finer covers.
Thus, we may assume that for each $n$, the cover $\cU_{n+1}$ refines the cover $\cU_n$.

Fix an injective sequence $\seq{x}$ in $X$.
For each $n$, let
\[
\cV_n :=\set{U\sm\{x_n\}}{U\in\cU_n}.
\]
Then $\cV:=\UnN{n}\cV_n\in\Om(X)$, and thus there is a point-cofinite subfamily $\cW\sub\cV$.
It follows that for each $n$, there are only finitely many members of the family $\cV_n$ in $\cW$.
By moving to a subset of $\cW$, we may assume that there is a set $I\in\roth$ and for each $n\in I$ a set $U_n\in\cU_n$ with
\[
\cW = \set{U_n\sm\{x_n\}}{n\in I}.
\]
For $n\notin I$, let $k\in I$ be greater than $n$, and choose a set $U_n\in\cU_n$ with $U_k\sub U_n$.
By Lemma~\ref{lem:om}(2), we have $\{\seq{U}\}\in\Ga(X)$.
\epf

In particular, the property $\omg$ implies $\sone(\Ga,\Ga)$;
we will use this connection.

\bdfn
A set $S\sub\roth$ is \emph{centered} if every finite subset of $S$ has an infinite intersection.
An infinite set $a$ is a \emph{pseudointersection} of a set $S\sub\roth$ if for each set $s\in S$ we have $a\subs s$.
Let $\fp$ be the minimal cardinality of a centered set that has no pseudointersection.
\edfn

It is known that $\aleph_1\le\fp\le \fb$, and that strict inequalities are consistent~\cite{BlassHBK}.

\blem[Arhangel'ski\u{\i}]
\label{lem:w-lind}
Every open $\w$-cover of a set of real numbers has a countable $\w$-subcover of that set.
\elem
\bpf
Let $\cU\in\Om(X)$.
For each natural number $k$, we have $\sset{U^k}{U\in\cU}\in\Op(X^k)$.
Since $X$ is a set of real numbers, its finite powers are Lindel\"of, and there is a countable subcover
$\cU_k\in\Op(X^k)$.
It follows that for each $k$-element set $F\sub X$ there is an element $U^k\in\cU_k$ with $F\sub U$.
Thus, we have
\[
\UnN{k}\sset{U}{U^k\in\cU_k}\in\Om(X).\qedhere
\]
\epf

\bprp[Gerlits--Nagy]
\label{prp:nonomg}
$\non(\omg)=\fp$.
\eprp
\bpf
Let $S\sub\roth$ be a centered set of cardinality $\fp$ that has no pseudointersection.
For each $n$, let
\[
O_n := \sset{a\in \PN}{n\in a}
\]
Since the set $S$ is centered, we have $\sset{O_n}{n\in\bbN}\in\Om(X)$.
If $\sset{O_n}{n\in b}\in\Ga(X)$, then the set $b$ is a pseudointersection of $S$.
A contradiction.
This shows that $\non(\omg)\le \fp$.

On the other hand, let $X$ be a set of real numbers with $\card{X}<\fp$.
Suppose that $\cU\in\Om(X)$.
By Lemma~\ref{lem:w-lind}, we may assume that the cover $\cU$ is countable.
Enumerate $\cU=\sset{U_n}{n\in\bbN}$.
The sets
\[
s_x := \sset{n\in\bbN}{x\in U_n},
\]
for $x\in X$, form a centered family.
Let $a$ be a pseudointersection of this family.
Then $\sset{U_n}{n\in a}\in\Ga(X)$.
\epf

\blem
\label{lem:ptower}
There is an unbounded $\fp$-tower if and only if $\fp=\fb$.
\elem
\bpf
$(\Impl)$ Let $T$ be an unbounded $\fp$-tower.
Since the set $T$ is unbounded, we have
\[
\fb\le\card{T}\le\fp\le\fb.
\]

$(\Leftarrow)$ Assume that $\fp=\fb$.
Then there is an unbounded set $\sset{f_\alpha}{\alpha<\fp}\sub\NN$.
We construct an unbounded tower $\sset{t_\alpha}{\alpha<\fp}\sub\roth$ by induction on $\alpha$.

Pick a set $t_0\in\roth$ with $f_0\le t_0$.
For $\alpha>0$, assume that we have defined $t_\beta$ for all $\beta<\alpha$, such that the sets $t_\beta$ are
$\subs$-decreasing with $\beta$.
Then the set $\sset{t_\beta}{\beta<\alpha}$ is centered, and its cardinality is smaller than $\fp$.
Let $t$ be a pseudointersection of the set $\sset{t_\beta}{\beta<\alpha}$, and pick an infinite set
$t_\alpha\sub t$ with $f_\alpha\le t_\alpha$.
\epf

\blem
\label{lem:sinf}
Let $X$ be an $\sogg$ space.
For each sequence $\seq{\cU}\in\Ga(X)$, there are infinite sets
$\cI_1\sub\cU_1$, $\cI_2\sub\cU_2$, \dots, such that
$\UnN{n}{\cI_n}\in\Ga(X)$.
\elem
\bpf
Split every cover $\cU_n$ into (countably) infinitely many infinite, disjoint subsets.
These subsets are all in $\Ga(X)$.
Arrange all of these sets (for all $n$) into a single sequence, and apply $\sogg$.
\epf

\blem
\label{lem:pd}
Let $A_1,A_2,\dotsc$ be infinite sets.
There are infinite pairwise disjoint subsets $B_1\sub A_1, B_2\sub A_2, \dotsc$.
\elem
\bpf
Begin with the sets $B_n:=\emptyset$ for all $n$.
For each $k=1,2,3,\dotsc$, for each $i=1,\dotsc,k$, pick an element $a\in A_i\sm\UnN{n} B_n$ and
add it to the set $B_i$.
\epf

The following result of Jordan is crucial for understanding $\omg$ sets~\cite[Theorem~7]{Jordan08}.

\bprp[Jordan]
\label{prp:Jordan}
Let $X=\UnN{n}X_n$ be an increasing union of $\sogg$ spaces.
Suppose that
\[
\cU_1\in\Gamma(X_1),\quad
\cU_2\in\Gamma(X_2),\quad
\cU_3\in\Gamma(X_3),\quad
\dotsc.
\]
Then there are sets $U_1\in\cU_1, U_2\in\cU_2,\dotsc$ with
$\sset{U_n}{n\in\bbN}\in\Ga(X)$.
\eprp
\bpf
By thinning out the given families, we may assume that they are pairwise disjoint (Lemma~\ref{lem:pd}).

Step 1: By Lemma~\ref{lem:sinf}, we may thin out the families $\seq{\cU}$, such that they remain infinite,
and $\UnN{n}\cU_n\in\Ga(X_1)$.
Fix a set $U_1\in\cU_1$.

Step 2: By the same Lemma, we may thin out further the families $\cU_2,\cU_3,\dotsc$, such that
they remain infinite, and $\Un_{n=2}^\oo\cU_n\in\Ga(X_2)$.
Fix a set $U_2\in\cU_2$.

Step $k$: By Lemma~\ref{lem:sinf}, we may thin out further the families $\cU_k,\cU_{k+1},\dotsc$, such that
they remain infinite, and $\Un_{n=k}^\oo\cU_n\in\Ga(X_k)$.
Fix a set $U_k\in\cU_k$.

For each $k$, we have
\[
\{\seq{U}\}\subs \Un_{n=k}^\oo\cU_n\in\Ga(X_k),
\]
and thus $\{\seq{U}\}\in\Ga(X_k)$ for all $k$.
It follows that $\{\seq{U}\}\in\Ga(X)$.
\epf

We are now ready to prove the most important
result regarding the existence of sets of real numbers with selective covering properties~\cite[Theorem~3.6]{LinSAdd}.
\footnote{While I firmly believe it, this assertion is subjective.
If it helps, the referee report reads ``I can't think of any more breakthrough
construction than~\cite[Theorem~3.6]{LinSAdd} in this area.''}

\bthm[Orenshtein--Tsaban]
\label{thm:crown}
For each unbounded $\fp$-tower $T\sub\roth$,
the set $T\cup\Fin$ is $\omg$.
\ethm
\bpf
Let $T = \sset{t_\alpha}{\alpha<\fp}$ be an unbounded $\fp$-tower.
Let $X := T\cup\Fin$.
For an ordinal $\alpha<\fb$, let
\begin{align*}
X_{<\alpha} & := \Fin\cup\sset{t_\beta}{\beta<\alpha};\\
X_{\ge\alpha} & := \Fin\cup\sset{t_\beta}{\alpha\le \beta<\fb}.
\end{align*}
We will use, repeatedly, the following simple variation of observations made in the earlier proofs.

\blem
\label{lem:ab}
Let $\cU\in\Om(T\cup\Fin)$.
For each ordinal $\alpha<\fp$, there are an ordinal $\beta<\fp$ with $\alpha<\beta$ and a set $\cV\sub\cU$ with
$\cV\in\Ga(X_{<\alpha}\cup X_{\ge\beta})$.
\elem
\bpf
If $\cV\in\Ga(X_{<\alpha}\cup X_{\ge\beta_0})$, then
$\cV\in\Ga(X_{<\alpha}\cup X_{\ge\beta})$ for all ordinals $\beta<\fp$ with $\beta_0\le\beta$.
Thus, it suffices to prove the lemma without the restriction $\alpha<\beta$.

Since $\card{X_{<\alpha}}<\fp$ and $\cU\in\Om(X_{<\alpha})$,
there is a set $\cU'\sub\cU$ with $\cU'\in\Ga(X_{<\alpha})$ (Lemma~\ref{prp:nonomg}).
By Lemma~\ref{lem:gm}, there are a set $a\in\roth$ with $a(1)=1$ and distinct sets $\seq{U}\in\cU'$ such that for each $n$,
all points $x\in\PN$ that omit the interval $\ooi{a}{n}$ are in the set $U_n$.
We have $\{\seq{U}\}\in\Ga(X_{<\alpha})$.

Since the set $T$ is unbounded, there is an ordinal $\beta<\fp$ such that the set
\[
I := \sset{n}{t_\beta\cap\ooi{a}{n}=\emptyset}
\]
is infinite (Lemma~\ref{lem:bddslalom}).
Then $\cV := \sset{U_n}{n\in I}\in \Ga(X_{<\alpha})$.

Let $\gamma<\fp$ be an ordinal with $\beta\le\gamma$.
For almost all $n\in I$, we have
\[
t_\gamma\cap\ooi{a}{n}\sub t_\beta\cap\ooi{a}{n}=\emptyset,
\]
and thus $t_\gamma\in U_n$. Thus, we also have $\cV \in \Ga(X_{\ge\beta})$.
\epf

Let $\alpha_1:=0$.
By Lemma~\ref{lem:ab}, for each natural number $n$ there are a set $\cV_n\sub\cU$ and an ordinal $\alpha_{n+1}<\fp$
such that $\alpha_n<\alpha_{n+1}$ and
$\cV_n\in\Ga(X_{\alpha_{n}}\cup X_{\ge\alpha_{n+1}})$.
Let $\alpha := \sup_n \alpha_n$.
For each $n$, let
\[
Y_n := X_{\alpha_{n}}\cup X_{\ge\alpha}.
\]
Then $X = \UnN{n}Y_n$, the union is increasing, and $\cV_n\in\Ga(X_n)$ for all $n$.

For each ordinal $\gamma\le\alpha$, the set
$\sset{\beta}{\beta<\gamma}\cup\sset{\beta}{\alpha\le\beta<\fb}$
is order-isomorphic to the set $\sset{\beta}{\beta<\fb}$, and thus
the set
\[
\sset{t_\beta}{\beta<\gamma}\cup\sset{t_\beta}{\alpha\le\beta<\fp}
\]
is an unbounded $\fb$-tower.
By Theorem~\ref{thm:btowersogg}, for each $n$, the set $Y_n$ is $\sogg$.

By Proposition~\ref{prp:Jordan}, there are sets
$U_1\in\cV_1$,
$U_2\in\cV_2$,
$U_3\in\cV_3$,
\dots,
such that
$\{\seq{U}\}\in\Ga(X)$.
\epf

\subsection{Comments for Section~\ref{sec:crown}}

The interest in the Crown Theorem (Theorem~\ref{thm:crown}) is as old as the Gerlits--Nagy seminal paper~\cite{GN82}.
The property $\omg$ is very restrictive.
We have
\[
\omg=\sone(\Om,\Ga) \longrightarrow \sone(\Om,\Op) = \sone(\Op,\Op).
\]
For the last equality~\cite[Theorem~17]{coc1}, suppose that our space is $\sone(\Om,\Op)$ and we are given a sequence of open covers.
We split the sequence into infinitely many disjoint subsequences, and produce from each subsequence a single $\omega$-cover,
consisting of all finite unions taking at most one set from each open cover.
We then apply our property $\sone(\Om,\Op)$ to obtain a cover, and disassemble the unions to single sets.
This shows that our space is $\sone(\Op,\Op)$.
In particular, $\omg$ sets of real numbers are strong measure zero.

Gerlits and Nagy noticed that uncountable $\omg$ sets of real numbers exist
if $\aleph_1<\fp$ (Proposition~\ref{prp:nonomg}), but did not know whether they may be nontrivial.
Galvin and Miller~\cite[Theorem~1]{GM84} proved that nontrivial $\omg$ sets exist if $\fp=\fc$.
Just, Miller, Scheepers and Szeptycki~\cite[Theorem~5.1, Case~2]{coc2} proved that for each unbounded $\aleph_1$-tower
$T$, the set $T\cup\Fin$ is $\sone(\Omega,\Omega)$, and Scheepers~\cite[Theorem~3]{wqn} proved that such a set must be $\sone(\Gamma,\Gamma)$.
They asked whether it must be $\omg$~\cite[Problem~7]{coc2}.
Scheepers proved that for each unbounded $\fp$-tower $T$, the set $T\cup\Fin$ is $\sone(\Gamma,\Gamma)$~\cite[Theorem~6]{Scheepers98}.
In Theorem~11 of the ArXiv version of his paper~\cite{MillerBC},
Miller proves a theorem of Bartoszy\'nski that a property strictly stronger than $\aleph_1=\fb$, denoted $\diamondsuit(\fb)$, implies that there is a nontrivial $\omg$ set,
and asked whether the hypothesis $\aleph_1 = \fb$ suffices.
In their study of a problem of Malykhin, Gruenhage and Szeptycki~\cite[Corollary~10]{FUfin} prove that if $\fp=\fb$ then there is a set of real numbers with a property weaker than
$\omg$,  and wonder whether $\fp = \fb$ implies the existence of a nontrivial $\omg$ set.
The Crown Theorem improves upon all of these partial results, and answers these questions.

The author still does not fully understand his first proof of the Crown Theorem~\cite[Theorem~3.6]{LinSAdd}.
In 2016, while preparing a course on omission of intervals, he came up with a simpler proof, that he
presented in several conferences and included in a paper with Osipov and Szewczak~\cite[Theorem~6]{SSS}.
That proof was still somewhat technical.
He arrived at the present, conceptual proof while working on the present paper,
motivated by the preparation of his 2023 course on omission of intervals.

The property $\omg$ is preserved by finite powers.
It is not provably preserved by finite products:
Miller, Zdomskyy and the author~\cite[Theorem~3.2]{gammaprods} proved that,
assuming \CH{}, there are two $\omg$ (even for countable Borel covers)
sets of real numbers whose product is not even $\sfoo$. On the other hand~\cite[Theorem~2.8]{gammaprods} they proved that
for each unbounded $\aleph_1$-tower $T$, the product of the set $T\cup\Fin$
with every $\omg$ set is $\omg$.
Szewczak and W\l{}udecka~\cite[Theorem~4.1]{SzWl21} prove that
this is, in fact, the case for every (generalized) unbounded $\fp$-tower $T$.

\section{Refined scales and nonimplications}
\label{sec:nonimpl}

\subsection{Refinements}

The proofs of our main theorems tell more than their statements.
For the following applications, we need to state some of these theorems more generally.

\bdfn
Let $R$ be a binary relation on the set $\roth$.
A set $B\sub\roth$ is an \emph{$R$-refinement} of a set $A\sub\roth$ if
there is a surjection $\varphi\colon A\to B$ with $aR\varphi(a)$ for all $a\in A$.
\edfn

Every $\supseteq^*$-refinement is a $\les$-refinement and a $\sqe$-refinement.

\bthm
\label{thm:uidggp}
For each $\les$-refinement $S'$ of a $\fb$-scale $S$, the set
$S'\cup \Fin$ is $\uidgg$.
\ethm
\bpf
Let $S = \sset{s_\alpha}{\alpha<\fb}\sub\roth$ be a $\fb$-scale.
Using the surjection, we may write
$S' = \sset{s'_\alpha}{\alpha<\fb}$, where
$s_\alpha \les s'_\alpha$ for all $\alpha<\fb$.

We define a set $I$ as in our proof of Theorem~\ref{thm:bs}.
In the remainder of the proof, we change everywhere $s$ to $s'$,
except for the case $\beta\ge\alpha$, where we use that
$s_\alpha \les s_\beta \les s'_\beta$, so that
$s_\alpha \les s'_\beta$.

This establishes $\ufog$.
The upgrade to $\uidgg$ is as in the proof of Theorem~\ref{thm:un}, that is, using $\sogg$ instead of $\ufog$ in the first step.
\epf

Similar modifications of the proofs of Theorem~\ref{thm:utgg}, Theorem~\ref{thm:btowersogg} and Theorem~\ref{thm:crown}
establish the following results.
Here too, we use the transitivity of the relations $\sqe$ and $\subs$.

\bthm
\label{thm:utggp}
For each $\sqe$-refinement $S'$ of a $\fb(\sqe)$-scale $S$, the set
$S'\cup \Fin$ is $\utgg$. \qed
\ethm

\bthm
\label{thm:btowersoggp}
For each $\supseteq^*$-refinement $T'$ of an unbounded $\fb$-tower $T$, the set
$T'\cup \Fin$ is $\sogg$.
And if $\fp=\fb$, then this set is $\omg$.
\qed
\ethm

\subsection{Hitting sets}

\bdfn
We fix a partition $\bbN = \UnN{n}I_n$ with
\[
1 = \min I_1 < \min I_2 < \min I_3 < \dotsb.
\]
Let
\[
\cI := \set{s\in \PN}{s\cap I_n\text{ is finite for all }n}
\]

For each $n$ and each $m\in I_n$, let
\[
O^n_m := \sset{x \in \PN}{n\notin x}.
\]
Then $\cO_n := \sset{O^n_m}{m\in I_n}\in\Ga(\cI)$ for all $n$.

For a set $A$ and a natural number $k$, let $[A]^k := \sset{F\sub A}{\card{F}=k}$.
Let $f\in\NN$.
A set $S\sub \cI$ is \emph{$f$-hitting} if
for each function $g\in\prod_{n=1}^\oo [I_n]^{f(n)}$ there is an element $s\in S$ with $g(n)\sub s$ for infinitely many $n$.
It is \emph{$k$-hitting} if it is $f$-hitting for the constant function $f(n)=k$.
For $k=1$, we identify the sets $[I_n]^1$ with the sets $I_n$.
A set $S\sub \cI$ is \emph{$*$-hitting} if for each function $g\in\prod_{n=1}^\oo \Fin(I_n)$ there is an element $s\in S$ with $g(n)\sub s$ for infinitely many $n$.
\edfn

\bprp
\label{prp:hitting}
Let $X\sub\cI$.
\be
\item If the set $X$ is $1$-hitting, then it is not $\sogg$.
\item If the set $X$ is $k$-hitting for $k\ge 2$, then it is not $\mathsf{U}_k(\Ga,\Ga)$.
\item If the set $X$ is $f$-hitting for an unbounded function $f\in\NN$, then it is not $\uidgg$.
\item If the set $X$ is $*$-hitting, then it is not $\ufog$.
\ee
\eprp
\bpf
(1) Suppose that there are $m_1\in I_1$, $m_2\in I_2$, \dots such that
\[
\{O^1_{m_1},O^2_{m_2},\dotsc\}\in\Ga(X).
\]
Let $g(n) := m_n $ for all $n$.
For each set $x\in X$, for almost all $n$ we have $x\in O^n_{m_n}$, and thus $g(n) \notin x$.
A contradiction.

(2) Suppose that there are sets $F_1\in [I_1]^k$, $F_2\in [I_2]^k$, \dots such that
\[
\smallmedset{\Un_{m\in F_n} O^n_m}{n\in\bbN}\in\Ga(X).
\]
Let $g(n) := F_n$ for all $n$.
For each set $x\in X$, for almost all $n$,
there is $m\in F_n$ with $x\in O^n_m$, that is, $m\notin x$, and thus $g(n) \nsubseteq x$.
A contradiction.

(3) By Lemma~\ref{lem:limsup}, we have $\mathsf{U}_f(\Ga,\Ga)=\uidgg$.
The proof that the set $X$ is not $\mathsf{U}_f(\Ga,\Ga)$ is similar to the previous proofs.

The proof (4) is also similar.
\epf

\subsection{Nonimplications}

\bthm
\label{thm:u2nots1}
Assume that $\fb=\fc$.
Every $\sqe$-refinement $S$ of a $\fb(\sqe)$-scale has a $\sqe$-refinement $S'$ such the set $S'\cup\Fin$ is $\utgg$ but not $\sogg$.

Moreover, if $\fp=\fc$ then the set $S'\cup\Fin$ is also $\sone(\Om,\Om)$.
\ethm
\bpf
Let $S=\sset{s_\alpha}{\alpha<\fb}\sub\roth$ be a $\sqe$-refinement of a $\fb(\sqe)$-scale.
Enumerate
$\prod_{n=1}^\oo I_n = \sset{g_\alpha}{\alpha<\fc}$.
For each $\alpha<\fb$, pick an infinite subset $s'_\alpha\sub g_\alpha[\bbN]$ with $s_\alpha\sqe s'_\alpha$.
Let
\[
X := \sset{s'_\alpha}{\alpha<\fb}\cup \Fin.
\]
By Theorem~\ref{thm:utggp}, the set $X$ is $\utgg$.

The set $X$ is 1-hitting: For each function $g\in\prod_{n=1}^\oo I_n$,
there is $\alpha<\fc$ with $g = g_\alpha$.
Since $s'_\alpha\sub g_\alpha[\bbN]$, we have $g(n)\in s'_\alpha$ for infinitely many $n$.
By Proposition~\ref{prp:hitting}, the set $X$ it is not $\sogg$.

Assume that $\fp=\fc$.
The result we establish in this case is stronger than the asserted one.
It involves the Gerlits--Nagy $\delta$ property~\cite[Section~2]{GN82},
a property stronger than $\sone(\Om,\Om)$~\cite[Lemma~4.4]{GNP}.
S. Bardyla, J. \v{S}upina and L. Zdomskyy proved that,
if $\fp=\fc$, then every set $A\sub\roth$ of cardinality $\fc$ has a $\supseteq^*$-refinement $A'$
such that the set
$A'\cup\Fin$ is $\delta$~\cite[Theorem~3.1]{BSZ23}.

Let $S''$ be a $\supseteq^*$-refinement of the set $S'$ such that the set $S''\cup\Fin$ is a $\delta$-set (and thus $\sone(\Om,\Om)$).
By Theorem~\ref{thm:utggp}, the set $S''\cup\Fin$ is $\utgg$.
Since there are in the set infinite subsets of all possible images $g_\alpha[\bbN]$,
the set remains 1-hitting, and thus not $\sogg$ (Proposition~\ref{prp:hitting}).
\epf

\blem
\label{lem:kun}
Let $k\in\bbN$, $g\in  \prod_{n=1}^\oo [I_n]^k$ and $s\in\roth$.
There is an infinite set $t\sub g[\bbN]$ with $s\le s' := \Un t$.
\elem
\bpf
Let $\seq{m}$ be an increasing sequence with
$s(kn) < \min I_{m_n}$ for all $n$.
Let $s' := \UnN{n} g(m_n)$.
For each $n$, we have
\[
\card{\cbbl 1,s(kn)\cbbr\cap s'} \le \card{ \Un_{i=1}^{n-1} g(m_i) } = k(n-1),
\]
and thus $s(kn) < s'(k(n-1)+1)$.
For each $i$, pick $n$ with $i\in \cbbl k(n-1)+1,kn\cbbr$.
Then
\[
s(i)\le s(kn)\le s'(k(n-1)+1)\le s'(i). \qedhere
\]
\epf

\bthm
\label{thm:uidnotuk}
Assume that $\fb=\fc$.
Every $\les$-refinement $S$ of a $\fb$-scale has a $\les$-refinement $S'$ such that
the set $S'\cup\Fin$ is $\uidgg$ but not $\mathsf{U}_k(\Ga,\Ga)$ for any natural number $k$.

Moreover, if $\fp=\fc$ then the set $S'\cup\Fin$ is also $\sone(\Om,\Om)$.
\ethm
\bpf
Let $S=\sset{s_\alpha}{\alpha<\fb}\sub\roth$ be a $\les$-refinement of a $\fb$-scale.
Enumerate
$\UnN{k}\prod_{n=1}^\oo [I_n]^k = \sset{g_\alpha}{\alpha<\fc}$.
For each $\alpha<\fb$, pick an infinite subset $t_\alpha\sub g_\alpha[\bbN]$ with $s_\alpha\les s'_\alpha := \Un t_\alpha$ (Lemma~\ref{lem:kun}).
Let
\[
X := \sset{s'_\alpha}{\alpha<\fb}\cup \Fin.
\]
By Theorem~\ref{thm:uidggp}, the set $X$ is $\uidgg$.
By the construction, this set is $k$-hitting, and thus not $\mathsf{U}_k(\Ga,\Ga)$, for all $k$ (Proposition~\ref{prp:hitting}).

If $\fp=\fc$, then we can thin out the sets $t_\alpha$ further so that the resulting set $X$ is a $\delta$ set~\cite[Proof of Theorem~4.2]{DeltaSets} and, in particular, $\sone(\Om,\Om)$.
The set $X$ remains $k$-hitting for all $k$, and still $\les$-refines a $\fb$-scale.
\epf

\subsection{Comments for Section~\ref{sec:nonimpl}}
Zdomskyy and the author~\cite[Theorem~2.1]{HurPerf} proved that if $\fb=\fc$
then there is a $\fb$-scale $S$ such that the set $S\cup\Fin$ is $\ufog$ but not $\sogg$.
The author~\cite[Theorem~4.5]{MHP} proved that, assuming $\fb=\fc$, there is a $\fb(\sqe)$-scale $S$ such that the set $S\cup\Fin$ is $\utgg$ but not $\sogg$.
Theorem~\ref{thm:u2nots1} improves and extends this result.
Its proof, using omission of intervals, is substantially simpler.

Liu, He and Zhang~\cite[Theorem~4.1]{LHZ22} proved that if $\fb=\fc$ then there is a set of real numbers that is $\uidgg$ but not $\mathsf{U}_k(\Ga,\Ga)$ for any $k$.
Their proof is technical and uses a 2-dimensional version of omission of intervals.
Theorem~\ref{thm:uidnotuk} improves and extends this result,
using a much simpler proof.

The result of Bardyla--\v{S}upina--Zdomskyy~\cite[Theorem~3.1]{BSZ23} employed in the proofs of Theorems~\ref{thm:u2nots1} and~\ref{thm:uidnotuk}
is included in their proof of their breakthrough result that $\delta$ sets need not be $\omg$.
This settled an important problem from the Gerlits--Nagy seminal paper~\cite{GN82}.
Bardyla, \v{S}upina and Zdomskyy proceed to establish an impressive array of results of this type,
classifying large classes of related properties, using omission of intervals among other methods.
The author has recently wrote a variation of the Bardyla--\v{S}upina--Zdomskyy proof that gives a stronger result~\cite{DeltaSets}.

The Crown Theorem asserts that the hypothesis $\fp=\fb$ suffices for the existence of a nontrivial $\omg$ set of real numbers.
The examples considered in this paper are all based on unions of (variations of) $\fb$-scales with the set $\Fin$.
The hypothesis $\fp=\fb$ is consistent with the assertion that every union of a $\fb$-scale with $\Fin$ is an $\omg$ set~\cite{CRZ15}.
It remains open whether counterexamples such as the ones presented in this section can be proved to exist based on the hypothesis $\fp=\fb$.
A closely related question is whether it is provable, outright in ZFC, that the properties $\omg$ and $\sogg$
are nonequivalent for subsets of the real line.

\section{Perturbations}

By our results, the types of sets considered in the previous sections cannot provide examples of the kinds discussed below.
For that, we consider perturbations of scales.

\subsection{When $k$ sets are not enough but $k+1$ are}
Replacing 2 by $k$ in the proof of Lemma~\ref{lem:gm2} establishes the following observation.

\blem
\label{lem:gmk}
Assume that $\Fin\sub X\sub \PN$.
For each cover $\cU\in\Om(X)$,
there are a function $b\in\roth$ with $b(1)=1$ and distinct sets $\seq{U}\in\cU$
such that, for each $n$,
all points $x\in\PN$ with
\[
\card{x\cap\ooi{b}{n}}\le k,
\]
are in the set $U_{(k+1)(n-1)+1}\cup \dotsb \cup U_{(k+1)n}$.\qed
\elem

\blem
\label{lem:2int}
Let $s\in\roth$, $k\in\bbN$ and $g\in \prod_{n=1}^\oo [I_n]^k$.
There is an infinite set $t\sub g[\bbN]$
such that the set $s' := \Un t$ has
\[
\card{s'\cap \obbl s(n),s(n+2) \obbr}\le k
\]
for all $n$.\qed
\elem

\bthm
\label{thm:uk+1notuk}
Assume that $\fb=\fc$.
For each $k$,
there are $\mathsf{U}_{k+1}(\Ga,\Ga)$ sets of real numbers that are not $\mathsf{U}_k(\Ga,\Ga)$.
\ethm
\bpf
Let $k\in\bbN$.
Enumerate $\prod_{n=1}^\oo [I_n]^k = \sset{ g_\alpha }{\alpha<\fc}$.
Let $S = \sset{s_\alpha}{\alpha<\fc}$ be a $\fb(\sqe)$-scale.
For each $\alpha<\fc$, pick an infinite set $t_\alpha\sub g_\alpha[\bbN]$
such that the set $s'_\alpha := \Un t_\alpha$ has
\[
\card{s'_\alpha\cap \obbl s_\alpha(n),s_\alpha(n+2) \obbr}\le k
\]
for all $n$ (Lemma~\ref{lem:2int}).
Let $X := \sset{ s'_\alpha}{\alpha<\fc}\cup\Fin$.

The set $X$ is $k$-hitting, and thus not $\mathsf{U}_k(\Ga,\Ga)$.
To see that the set is $\mathsf{U}_{k+1}(\Ga,\Ga)$, we proceed as in the proof of Theorem~\ref{thm:utgg}, using Lemma~\ref{lem:gmk}.

Let $\seq{\cU}\in\Ga(S\cup\Fin)$.
By Lemma~\ref{lem:gmk},
for each $n$ there are a function $b_n\in\roth$ with $b_n(1)=1$ and distinct sets $\seq{U^n}\in\cU_n$
such that, for each $m$,
all points $x\in\PN$ with
\[
\card{x\cap\ooi{b_n}{m}}\le k,
\]
are in the set
\[
U^n_{(k+1)(m-1)+1}\cup \dotsb \cup  U^n_{(k+1)m}.
\]
By Lemma~\ref{lem:bdd2}, there is a set $b\in\roth$ with $\sset{b_n}{n\in\bbN}\sqe b$.
There is a function $c\in\roth$ such that, for each $n$, for all $m\ge c(n)$,
the interval $\ooi{b}{m}$ contains an open $b_n$-interval,
say $\ooi{b_n}{i_{n,m}}$.
If $m\ge c(n)$ and a set $x\in\PN$ has
\[
\card{x\cap\ooi{b}{m}}\le k,
\]
then $\card{x\cap\ooi{b_n}{i_{n,m}}}\le k$,
and thus
\[
x\in W^n_m := U^n_{(k+1)(i_{n,m}-1)+1}\cup \dotsb \cup U^n_{(k+1)i_{n,m}}.
\]
As the set $S$ is unbounded, there is by Lemma~\ref{lem:bddslalom} an ordinal $\alpha<\fb$,
such that the set
\[
I := \sset{m}{s_\alpha \cap \ooi{b}{m} = \emptyset}
\]
is infinite.
Since the set $I$ is infinite, for each $n$ we have
\[
\set{W^n_{m}}{c(n)\le m\in I}\in\Ga(S\cup\Fin).
\]
The set $\Fin\cup\set{s_\beta}{\beta<\alpha}$ is $\sogg$ (Lemma~\ref{lem:nons1gg}).
Thus, there is a function $g\in\NN$ with $c(n)\le g(n)\in I$ for all $n$, such that
\[
\set{W^n_{g(n)}}{n\in\bbN}\in\Ga(\Fin\cup\set{s_\beta}{\beta<\alpha}).
\]
We have $n\le c(n)\le g(n)$ for all $n$.

Let $\beta\ge \alpha$. Then $s_\alpha\sqe s_\beta$.
As in the proof of Theorem~\ref{thm:utgg},
for almost all $n$ we have $\card{s_\beta\cap\ooi{b}{g(n)}}\le 1$.
Let $h\in\NN$ be a function with
\[
\ooi{b}{g(n)} \sub \obbl s_\beta(h(n)),s_\beta(h(n)+2)\obbr
\]
for almost all $n$.
Then for almost all $n$, we have
\[
\card{s'_\beta\cap\ooi{b}{g(n)}} \le \card{s'_\beta\cap \obbl s_\beta(h(n)),s_\beta(h(n)+2) \obbr} \le k
\]
and thus $s'_\beta\in W^n_{g(n)}$.
\epf

\subsection{When boundedly many sets are not enough but finitely many are}

We will use bounded perturbations of $\fb$-scales.

\blem
\label{lem:perturbedscale}
Let $g\in\prod_{n=1}^\oo [I_n]^{f(n)}$, where $f\in\NN$.
For each set $a\in\roth$ there is a set $c\in\roth$ such that:
For each set $s\in\roth$,
every $c$-interval omitted by the set $s$ contains an $a$-interval omitted by the set
\[
b := \Un g[s] = \Un_{n\in s}g(n).
\]
\elem
\bpf
By induction on $n$, choose a function $c\in\roth$ such that the interval $\ooi{c}{n}$
contains more than
\[
d(n) := \sum_{i=1}^{c(n)} f(i)
\]
open $a$-intervals.

Suppose that the set $s$ omits a $c$-interval $\ooi{c}{n}$.
Then
\[
s\cap \cbbl 1,c(n+1)\obbr\sub \cbbl 1,c(n)\obbr.
\]
For $i\ge c(n+1)$, since $g(i) \sub I_i$, we have
\[
c(n+1)\le i\le \min g(i).
\]
It follows that
\[
\card{b\cap \cbbl 1,c(n+1)\obbr}
\le
\card{\Un_{i=1}^{c(n)-1}g(i)}
=
\sum_{i=1}^{c(n)-1}f(i)
\le
d(n).
\]
Since the interval $\ooi{c}{n}$ contains more than $d(n)$ open $a$-intervals,
the set $b$ omits one of these $a$-intervals.
\epf

\bprp
\label{prp:kunb}
Let $S\sub\roth$ be a $\kappa$-unbounded set, for a cardinal $\kappa$ of uncountable cofinality.
For each element $s\in S$ let $g_s\in\prod_{n=1}^\oo [I_n]^{f(n)}$.
Then the set $S' := \sset{\Un g_s[s]}{s\in S}\sub\roth$ is $\kappa$-unbounded.
\eprp
\bpf
Since the sets $I_n$ are disjoint, the map $s\mapsto \Un g_s[s]$ is bijective.
In particular, we have $\kappa\le\card{S}=\card{S'}$.

Let $a\in\roth$.
Apply Lemma~\ref{lem:perturbedscale} to the set $\tilde a$ to obtain a set $c$.
Let $d(n) := c(2n)$ for all $n$.

Let $s\in S$.
For each $n$, if $c(2n)\le s(n)$, then the set $s$ omits at least $n$ $c$-intervals.
Assume that $d\lei s$.
Then $s$ omits infinitely many $c$-intervals, and thus the set $\Un g_s[s]$ omits infinitely many $\tilde a$-intervals.
It follows that $a\lei \Un g_s[s]$ (Corollary~\ref{cor:omit1}).
Thus,
\[
\sset{s\in S}{\Un g_s[s]\les a}
\sub
\sset{s\in S}{s\les d}.
\]
It follows that only fewer than $\kappa$ members of the set $S'$ are bounded by $a$
(Remark~\ref{rem:unccof}).
\epf

\bthm
\label{thm:ufognonuidgg}
Assume that $\fb=\fc$.
There is a $\fb$-unbounded set $S\sub\roth$ such that the set
$S\cup\Fin$ is $\ufog$ and $\sone(\Ga,\Op)$ but not $\uidgg$.
\ethm
\bpf
Let $\sset{s_\alpha}{\alpha<\fc}$ be a $\fb$-scale.
Enumerate $\prod_{n=1}^\oo [I_n]^n = \sset{g_\alpha}{\alpha<\fc}$.
We will verify that the set
\[
S := \sset{\Un g_\alpha[s_\alpha]}{\alpha<\fc}
\]
is as required.

Since $\fb$-scales are $\fb$-unbounded (Lemma~\ref{lem:bsc-bu}),
and the cardinal $\fb$ has uncountable cofinality (indeed, it is regular),
the set $S$ is $\fb$-ubounded (Proposition~\ref{prp:kunb}).
It follows that the set $S\cup\Fin$ is $\sone(\Ga,\Op)$ (Lemma~\ref{lem:dconc2}).

The proof that the set $S\cup\Fin$ is $\ufog$ is similar to our proof of Theorem~\ref{thm:bs}:
It suffices to prove that our set is $\ufin(\Ga,\Ga)$.
Let $\seq\cU\in\Ga(X)$.
For each $n$, take a function $a_n\in\roth$ and distinct sets $\seq{U^n}$ for the cover $\cU_n$ as in Lemma~\ref{lem:gm}.
For each $n$, apply Lemma~\ref{lem:perturbedscale} to define a set $c_n$ corresponding to the set $a_n$.

Pick an ordinal $\alpha<\fb$ be such that the set
\[
I:=\sset{n}{c_n(n+1) \le s_\alpha(n)}
\]
is infinite.
By Lemma~\ref{lem:nons1gg}, the set $\Fin\cup\sset{\Un g_\beta[s_\beta]}{\beta<\alpha}$ is $\sogg$.
Thus, there are sets $U^n_{m_n}\in\cU_n$, for $n\in I$, such that
\[
\sset{U^n_{m_n}}{n\in I}\in\Ga(\sset{\Un g_\beta[s_\beta]}{\beta<\alpha}).
\]
Let $\cF_n:=\emptyset$ for $n\nin I$, and
\[
\cF_n:=\sset{U^n_m}{\ooi{a_n}{m}\sub \cbbl 1,c_n(n+1)\obbr}\cup\{U^n_{m_n}\}
\]
for $n\in I$.
Then
\[
\sset{\Un\cF_n}{n\in\bbN}=\sset{\Un\cF_n}{n\in I}\cup\{\emptyset\}.
\]
By Lemma~\ref{lem:om}(2), it suffices to show that for each point $x\in X$, we have $x\in\Un\cF_n$ for almost all $n\in I$.

If $x\in \Fin\cup\sset{\Un g_\beta[s_\beta]}{\beta<\alpha}$, this is evident.

Fix an ordinal $\beta\ge\alpha$.
Since $s_\alpha \les  s_\beta$,
for almost all $n\in I$ we have:
\[
c_n(n+1)\le s_\alpha(n) \le s_\beta(n),
\]
and thus the set $s_\beta$ omits a $c_n$-subinterval of the interval $\cbbl 1,c_n(n+1)\obbr$.
By Lemma~\ref{lem:perturbedscale}, there is an $a_n$-subinterval
$\ooi{a_n}{m}$ of this $c_n$-interval
that is omitted by the set $\Un g_\beta[s_\beta]$.
We have $\ooi{a_n}{m}\sub \cbbl 1,c_n(n+1)\obbr$, and thus
$\Un g_\beta[s_\beta]\in U^n_m\in\cF_n$.

The set $S\cup\Fin$ is $\op{id}$-hitting:
Let $g\in\prod_{n=1}^\oo [I_n]^n$.
There is an ordinal $\alpha<\fc$ with $g_\alpha=g$.
Then $g(n)\sub \Un g_\alpha[s_\alpha]$ for all $n\in s_\alpha$.
By Proposition~\ref{prp:hitting}, the set $S\cup\Fin$ is not $\uidgg$.
\epf

\subsection{Comments for Section~\ref{sec:nonimpl}}

Liu, He and Zhang~\cite[Theorem~5.2]{LHZ22} established Theorem~\ref{thm:uk+1notuk} assuming \CH{}.
The present proof uses the weaker hypothesis $\fb=\fc$.
Similarly, Theorem~\ref{thm:ufognonuidgg} was earlier established by Liu, He and Zhang~\cite[Theorem~3.1]{LHZ22} assuming \CH{}.
In both cases, we obtain simpler proofs of stronger theorems.

\section{Concluding remarks}

In the author's view, \emph{omission of intervals} became a method when
\ifauthor
he
\else
it was
\fi
realized that it is not only relevant to $\omg$ sets~\cite{GM84} and the closely related property $\sogg$~\cite{coc2, Scheepers98}, but also to the Hurewicz property $\ufog$ and
to questions related to rate of growth~\cite{MHP}, such as the solution of the Hurewicz Problem presented here.
These were previously addressed by other, more topological means~\cite{BaSh01, ChaPol, sfh}.
The usage of omission of intervals led to the complete
analysis of the finer covering properties $\mathsf{U}_f(\Ga,\Ga)$,
as detailed here; we could not even dream about these achievements
before the development of this method.

The introduction of omission of intervals suggests a paradigm change:
Rather than asking about the \emph{existence} of sets of real numbers with certain covering properties, we often
prove that \emph{all sets} with an appropriate combinatorial structure have these properties.
The presented conceptual proof of the Crown Theorem, using that all intermediate sets are $\sogg$, and the simplicity and strength of the proofs in the last two sections, demonstrate the importance of this paradigm change.

This paper is the first to provide a solid foundation for omission of intervals.
The readers are now equipped with the needed tools to read any paper that employs this method, and to
continue this journey on their own.

\ed